\documentclass[a4paper, 10pt, oneside, reqno]{amsart}

\usepackage[utf8]{inputenc}
\usepackage{amsfonts,amssymb,amsthm,amsmath}
\usepackage{mathrsfs}
\usepackage{geometry}
\usepackage{graphicx,graphics}
\usepackage{color}
\usepackage{tikz-cd}
\usepackage[all]{xy}
\usepackage{comment}
\usepackage{subfiles}
\usepackage{stmaryrd}
\usepackage{hyperref}
\usepackage{todonotes}

\hypersetup{colorlinks = true, linkbordercolor = {white}, linkcolor = {blue}, citecolor = {blue}}

\theoremstyle{plain}
\newtheorem{thm}{Theorem}[section]
\newtheorem{lem}[thm]{Lemma}
\newtheorem{prop}[thm]{Proposition}
\newtheorem{corol}[thm]{Corollary}

\theoremstyle{definition}
\newtheorem{defi}[thm]{Definition}

\newtheorem{constr}[thm]{Construction}
\theoremstyle{remark}
\newtheorem{rque}[thm]{Remark}
\newtheorem{rques}[thm]{Remarks}

\newcommand{\F}{\mathbb{F}}
\newcommand{\Fq}{\F_q}

\newcommand{\Fqb}{\overline{\F}_q}
\newcommand{\mcO}{\mathcal{O}}
\newcommand{\Ql}{\mathbb{Q}_{\ell}}

\newcommand{\Z}{\mathbb{Z}}
\newcommand{\Spec}{\mathrm{Spec}}
\newcommand{\N}{\mathbb{N}}
\newcommand{\Frob}{\mathrm{F}}

\newcommand{\can}{\mathrm{can}}
\newcommand{\id}{\mathrm{id}}
\newcommand{\coev}{\mathrm{coev}}
\newcommand{\ev}{\mathrm{ev}}
\newcommand{\indlis}{\mathrm{indlis}}
\newcommand{\indcons}{\mathrm{indcons}}
\newcommand{\ULA}{\mathrm{ULA}}
\newcommand{\Perv}{\mathrm{Perv}}
\newcommand{\mcL}{\mathcal{L}}
\newcommand{\pr}{\mathrm{pr}}
\newcommand{\univ}{\mathrm{univ}}

\newcommand{\Cht}{\mathrm{Sht}}
\newcommand{\ChtR}{\Cht\mathcal{R}}
\newcommand{\Gr}{\mathrm{Gr}}
\newcommand{\mcS}{\mathcal{S}}
\newcommand{\mcG}{\mathcal{G}}
\newcommand{\Res}{\mathcal{R}}
\newcommand{\mcA}{\mathcal{A}}
\newcommand{\mcY}{\mathcal{Y}}
\newcommand{\mcF}{\mathcal{F}}
\newcommand{\mcE}{\mathcal{E}}
\newcommand{\mcB}{\mathcal{B}}

\newcommand{\DD}{\mathrm{D}}
\newcommand{\Rep}{\mathrm{Rep}}

\title{Nearby cycles commute with proper direct image on stacks of shtukas}
\author{Arnaud Eteve and Cong Xue}

\begin{document}

\maketitle

\noindent\textbf{Abstract.} Let $G$ be a generically reductive group over a smooth projective curve $X$ over a finite field. For any finite set $I$, we show that 
nearby cycles commute with proper direct image from stacks of shtukas to $X^I$. This generalizes some results of Salmon and the authors.

\tableofcontents

\section{Introduction}
In the introduction we illustrate the main result, under some simplification assumptions, for the stack of shtukas with one leg. In the next section we state the main theorems in the general setting, for the stacks of shtukas with several legs.

The stacks of shtukas with arbitrary many legs we considered in this paper are defined in \cite{Varshavsky} and recalled in \cite{Lafforgue}, which generalize the stacks of Drinfeld's shtukas with two legs.

\subsection{A simplified setting}
Let $X_0$ be a smooth projective curve over a finite field $\Fq$ of characteristic $p>0$. Let $x_0$ be a closed point (supposed of degree $1$) of $X_0$ and $N_0=nx_0$ a finite subscheme of $X_0$ for some $n \in \mathbb{N}$. 
Let $G$ be a connected split reductive group over $\Fq$. Let $^LG$ be the Langlands dual group of $G$ over $\mathbb{Q}_{\ell}$, where $\ell \neq p$. 

We denote by $X$ and $N$ the base change of $X_0$ and $N_0$ to $\Fqb$. We denote by $x \in X$ a geometric point over $x_0$.

Let $W$ be a representation of $^LG$. We denote by $\Cht_{\{1\}, W}$ the moduli stack of $G$-shtukas with one leg, which classifies $G$-bundles over $X$ and a modification (bounded by $W$) between the $G$-bundle and its inverse image of Frobenius. It has a projection to $X$ (morphism of leg). It is a Deligne-Mumford stack locally of finite type. Its local model is the Beilinson-Drinfeld affine Grassmannian $\Gr_{\{1\}, W}$.

We denote by $\Cht_{\{1\}, W, N}$ the moduli stack of $G$-shtukas with level structure $N$. It has a projection to $X-x$. 
We have the following commutative diagram:
$$
\xymatrix{
\Cht_{\{1\}, W, N} \ar[d]^{\pi} \\
\Cht_{\{1\}, W}  | _{(X-x)}   \ar@{^{(}->}[r]  \ar[d]^{\mathfrak{p}}
& \Cht_{\{1\}, W} \ar[d]^{\mathfrak{p}} \\
X-x  \ar@{^{(}->}[r]
& X 
}
$$
We are interested in the sheaf $$\mathcal{F} := \mathcal{S}_{\{1\}, W} \otimes \mathcal{L}^{nx}$$
defined over $\Cht_{\{1\}, W}  | _{(X-x)}$, 
where $\mathcal{S}_{\{1\}, W}$ is the Satake sheaf coming from the Beilinson-Drinfeld affine Grassmannian,
%defined over $\Cht_{\{1\}, W}$ 
and $\mathcal{L}^{nx}:=\pi_! \mathbb{Q}_{\ell}$ indicates
the level structure.
%is defined only over $\Cht_{\{1\}, W}  | _{(X-x)}$.
The geometric generic fiber of $\mathfrak{p}_! \mathcal{F}$ is the cohomology group with compact support of $\Cht_{\{1\}, W, N}$, which generalizes the space of automorphic forms with level $N$.

In \cite{XueSmoothness} we proved that the cohomology sheaf $\mathcal{H}_{G, N, \{1\}, W}:=\mathfrak{p}_! \mathcal{F}$ is ind-lisse over $X-x$.

In this paper we want to understand what happens at $x$. For this we consider the nearby cycle of $\mcF$ as follows.

\subsection{Nearby cycles for a trait}  \label{subsection-nearby-cycle-one-leg}
Let $X_{(x)}$ be the strict henselization of $X$ on $x$. Let $\overline{\eta}$ be a geometric generic point of $X$.
For any specialization map $\mathfrak{sp}: \overline{\eta} \rightarrow x$, i.e. a morphism $\overline{\eta} \rightarrow X_{(x)}$, we have a commutative diagram
$$
\xymatrix{
(\Cht_{\{1\}, W})|_{\overline{\eta}}  \ar[r]^{j}  \ar[d]^{\mathfrak{p}}
& \Cht_{\{1\}, W}  \ar[d]^{\mathfrak{p}}
& (\Cht_{\{1\}, W})|_{x}   \ar[l]_{i}  \ar[d]^{\mathfrak{p}} \\
\overline{\eta} \ar[r]
& X_{(x)}
& x   \ar[l]
}
$$
Note that $X_{(x)}$ is a trait. The classical nearby cycle functor for a trait (\cite[7, Exposé XIII]{SGA7}) attached to the specialization map $\mathfrak{sp}$
is defined to be
$$\Psi:=i^* j_*.$$ It is usually denoted by $R\Psi$.
Recall that we have a canonical morphism of functors 
\begin{equation}   \label{equation-p-Psi-to-Psi-p-one-leg}
\mathrm{can}: \mathfrak{p}_! \Psi  \rightarrow  \Psi \mathfrak{p}_! .
\end{equation}
coming from base change $\mathfrak{p}_! i^* = i^* \mathfrak{p}_!$ and $\mathfrak{p}_!  j_* \xrightarrow{\mathrm{unit}} j_* j^* \mathfrak{p}_!  j_* \simeq j_* \mathfrak{p}_! j^* j_* \xrightarrow{\mathrm{counit}} j_* \mathfrak{p}_!$ (where the middle isomorphism is base change).

%When $\mathfrak{p}$ is proper, (\ref{equation-p-Psi-to-Psi-p-one-leg}) is an isomorphism.
%When $\mathfrak{p}$ is not proper, (\ref{equation-p-Psi-to-Psi-p-one-leg}) may not be an isomorphism, i.e. the nearby cycle functor $\Psi$ may not commute with the proper direct image. 

%In general, the stack of shtukas $\Cht_{\{1\}, W}$ is not proper over $X$ (it is not even compact).
%However, we prove that 

\subsection{Main result (easiest case)}
 
\begin{thm}   \label{main-thm-one-leg}
For $\mcF$ defined as above, the canonical morphism induced by (\ref{equation-p-Psi-to-Psi-p-one-leg})
\begin{equation}   \label{equation-p-Psi-F-to-Psi-p-F-one-leg}
\mathrm{can}: \mathfrak{p}_! \Psi \mcF \rightarrow  \Psi \mathfrak{p}_! \mcF
\end{equation}
is an isomorphism.
\end{thm}

We denote by $\Cht = \Cht_{\{1\}, W}$. Theorem \ref{main-thm-one-leg} is equivalent to say that the canonical morphism 
\begin{equation}    \label{equation-H-c-Cht-x-to-H-c-Cht-eta} 
\mathrm{can}: H_c(\Cht | _{x}, \Psi \mathcal{F}) \rightarrow H_c(\Cht |_{\overline{\eta}}, \mathcal{F})
\end{equation}
is an isomorphism.

\subsection{Remarks}

\begin{rque}
Such kind of statement as in Theorem \ref{main-thm-one-leg} is false for general non proper schemes or stacks and general sheaves, i.e. in general the nearby cycles does not commute with the proper direct image.
\end{rque}

\begin{rque}
In most case, the stack of shtukas $\Cht_{\{1\}, W}$ is not proper over $X$ (it is not even compact).
But surprisingly Theorem \ref{main-thm-one-leg} is true for the sheaf $\mcF$. The reason relies on the properties of Satake sheaves: ULA (universally locally acyclic), functoriality, fusion; 
the partial Frobenius morphisms on stacks of shtukas; the twisted product structure of stacks of restricted shtukas.

The proof uses all these properties and a ``Zorro lemma" argument.
Note that the proof does not use any compactification.
\end{rque}

\begin{rque}
Theorem \ref{main-thm-one-leg} generalizes \cite{XueSmoothness}. And such type of statements appeared in \cite{Salmon}, \cite{Salmon2} and in the thesis of the first author \cite{EteveThese}. We send to \ref{subsection-lien-avec-Salmon} for a detailed discussion.
\end{rque}

\begin{rque}
In particular, Theorem \ref{main-thm-one-leg} implies that for any compactification of $\Cht$, the cohomology of the boundary with coefficient in $\Psi \mcF$ is zero: in fact, let $\overline{\Cht}$ be a compactification of $\Cht$ and $\iota: \Cht \rightarrow \overline{\Cht}$ the inclusion. We still denote by $\mcF$ the extension by zero to $\overline{\Cht}$. Then the open immersion $\iota$ induces a canonical morphism 
\begin{equation}   \label{equation-morphism-induced-by-open-immersion}
H_c(\Cht | _{x}, \Psi \mathcal{F} ) \rightarrow H(\overline{\Cht} | _{x}, \Psi \mathcal{F} )
\end{equation}
We can verify that (\ref{equation-morphism-induced-by-open-immersion}) coincides with (\ref{equation-H-c-Cht-x-to-H-c-Cht-eta}) (for the RHS, use $H(\overline{\Cht} | _{x}, \Psi \mathcal{F} ) \simeq H(\overline{\Cht} | _{\overline{\eta}}, \mathcal{F}) = H_c(\Cht | _{\overline{\eta}}, \mathcal{F})$). 
Note that the cone of (\ref{equation-morphism-induced-by-open-immersion}) is $H(\partial (\Cht)| _{x}, \Psi \mathcal{F} )$, where $\partial (\Cht) = \overline{\Cht} - \Cht$ is the boundary. Theorem \ref{main-thm-one-leg} is equivalent to say that the cone $H(\partial (\Cht)| _{x}, \Psi \mathcal{F} ) = 0$.
\end{rque}

\begin{rque}
Note that morphism (\ref{equation-H-c-Cht-x-to-H-c-Cht-eta}) is equivariant under the local Galois action, for the action of 
%$\mathrm{Gal}(\overline{F} / F)$ 
the global Galois group (of the function field of $X$)
on the cohomology group on the RHS and the action of 
%$\mathrm{Gal}(\overline{F}_x / F_x)$ 
the local Galois group (at the place $x$)
on the sheaf $\Psi \mathcal{F}$ on the LHS. 
%where $F$ is the function field of $X$ and $\overline{F}$ is the algebraic closure corresponding to $\overline{\eta}$, $F_x$ is the localization of $F$ at $x$ and $\overline{F}$ is an algebraic closure containing $\overline{F}$.

%By \cite{XueSmoothness}, the action of $\mathrm{Gal}(\overline{F} / F)$ on the cohomology group is unramified outside $x$. By Theorem \ref{main-thm-one-leg}, now we can study the ramification at $x$ by the action of $\mathrm{Gal}(\overline{F}_x / F_x)$ on $\Psi \mathcal{F}$.
\end{rque}

\subsection{Structure of this paper}
Main theorems of this paper Theorem \ref{thmMain2} and Theorem \ref{thmMain1} generalize Theorem \ref{main-thm-one-leg} to stacks of shtukas with several legs and more general sheaf $\mcF$, where $\mathcal{L}^{nx}$ is generalized by an arbitrary sheaf over a stack of restricted shtukas.

One reason to consider these more general sheaves is that we expect Theorem \ref{thmMain2} to be the main input in showing a strong form of local-global compatibility of the Langlands conjecture over function fields. This is one motivation of this paper. We send to \ref{subsection-local-global-compatibility} for more details.

We will give a proof of Theorem \ref{main-thm-one-leg} in section \ref{sectionOneLeg} to illustrate the general case. The general case Theorem \ref{thmMain2} and Theorem \ref{thmMain1} will be proved in sections \ref{sectionFusion}-\ref{sectionGeneralCase}.

\subsection{Acknowledgments}

The authors thank Vincent Lafforgue, Gérard Laumon, Alain Genestier, Dennis Gaitsgory and Jean-Fran\c cois Dat for many conversations and continuous support on this project. The first author was supported by the Max Planck Institute for Mathematics during the preparation of this paper.

\section{Notations and main theorems}

\subsection{Reminder on shtukas and sheaves}

\subsubsection{}
Let $X_0$ be a smooth projective curve over a finite field $\Fq$ of characteristic $p > 0$. Let $G_{\eta_0}$ be a connected reductive group over the generic point $\eta_0$ of $X_0$ and let $G$ be a smooth group scheme with geometrically connected fibers over $X_0$ with generic fiber $G_{\eta_0}$ and parahoric reduction at all the non reductive places. Let $R_0$ be the set of non reductive places. 

Let $\ell \neq p$ be a prime number, $E/\Ql$ a finite extension containing a square root of $q$ with ring of integers $\mathcal{O}_E$ and residual field $k_E$. Let $\Lambda \in \{E, \mathcal{O}_E, k_E\}$. Let us denote by ${^L}G$ the $L$-group of $G$. 

Let $N_0 \subset X_0$ be an effective divisor (i.e. a finite subscheme) and $x_0 \in N_0$. Let us denote by $X, N, R, ...$ the base change of all these $\Fq$-schemes to $\Fqb$. Let $\eta \in X$ be the generic point. Let $\bar{\eta}$ be a geometric point over $\eta$. We fix $\mathfrak{sp} : \bar{\eta} \to x$ a specialization map where $x$ is a geometric point of $X$ lies above $x_0$. Finally we denote by $\check{X} = X - (N \cup R)$ which is an open of $X$. 

\subsubsection{}
Let $I$ be a finite set, we have a stack of $G$-shtukas with $I$-legs $\Cht_{I} \xrightarrow{\mathfrak{p}} X^I$ as defined in \cite{Varshavsky} and
\cite{Lafforgue}. Above $\check{X}^I$ we have a finite étale cover $\Cht_{I,N} \xrightarrow{q_N} \Cht_I$. Let us denote by $\mcL_N = q_{N,*} \Lambda_{\Cht_{I,N}}$ the pushforward of the constant sheaf, this is a locally constant sheaf on $(\Cht_I)_{\check{X}^I}$. 

We denote by $\Gr_I = \Gr_{I,G}$ the Beilinson-Drinfeld Grassmannian over $X^I$ and by $L^+_IG\backslash \Gr_I$ its quotient by the global positive loop group. This quotient is called the Hecke-groupoid or the local Hecke stack. We furthermore have a local model map 
$$\varepsilon : \Cht_I \to L^+_IG\backslash \Gr_I$$
which is formally smooth \cite[Proposition 2.8]{Lafforgue}. 

The (ramified) geometric Satake correspondence \cite{ZhuRamifiedSatake}, \cite{RicharzGeoSatake} (original equivalence in \cite{MV} and the version for several legs in \cite{deJongConj}, the properties that we will need are recalled in \cite[Théorème 1.17, Théorème 12.16]{Lafforgue}), provides a collection of functors for all finite sets $I$
$$\mcS_I : \Rep_{\Lambda} {^L}G^I \to \Perv^{\ULA}((L^+_IG\backslash \Gr_I)_{(X-R)^I}, \Lambda)$$ 
where ${^L}G^I$ denote $I$-copies of ${^L}G$, and $\ULA$ denote the category of sheaves that are ULA (universally locally acyclic) relative to $(X-R)^I$. For any $W \in \Rep_{\Lambda} {^L}G^I$, we will still denote by $\mcS_{I,W}$ the pullback of the corresponding Satake sheaf on $(\Cht_I)_{|(X-R)^I}$. We will denote by $\Cht_{I,W} \subset \Cht_I$ the support of the sheaf $\mcS_{I,W}$. 

For a stack $\mcY$, we will denote by $\DD^b_c(\mcY, \Lambda)$ its usual bounded derived $\infty$-category of constructible sheaves and by $\DD_{\indlis}(\mcY, \Lambda) \subset \DD_{\indcons}(\mcY, \Lambda)$ its $\infty$-category of ind-lisse sheaves and ind-constructible sheaves in the sense of \cite{HemoRicharzScholbach}. In general we will call objects of these categories sheaves on $\mcY$. For a morphism $f : \mcY_1 \to \mcY_2$, whenever the functors are defined, we will denote by $f_!,f^!,f_*$ and $f^*$ the usual derived functors.

Finally the stack $\Cht_{I}$ is equipped with a filtration $\Cht_{I}^{\leq \mu}$ coming from Harder-Narasimhan trunctions where $\mu \in \Z^n$  \cite[Section 12]{Lafforgue}. Denote by $\mathfrak{p}^{\leq \mu} : \Cht_{I}^{\leq \mu} \to X^I$ the induced map. We denote by 
$$\mathfrak{p}_!  = \varinjlim_{\mu} \mathfrak{p}_!^{\leq \mu}$$ 
the functor from $\DD_c^b(\Cht_I, \Lambda)$ to $\DD_{\mathrm{indcons}}(X^I, \Lambda)$
%where the colimit is calculated in $\DD_{\mathrm{indcons}}(X^I, \Lambda)
(or to the category of sheaves on $\check{X}^I$ if this is relevant). 

\subsubsection{}  \label{subsection-def-F-I-W}

Let $W_0 \in \Rep_{\Lambda} {^L}G$ and $n \geq 0$, for $m \in \N$ large enough, there is a stack of restricted shtukas $\ChtR_{\{0\}, W_0,m, 0}^{nx}$ above $X$, as defined in \cite{GenestierLafforgue}. Let $(I,W)$ be a pair composed of finite set $I$ and $W \in \Rep_{\Lambda} {^L}G^I$. Let $(J,V)$ be another such pair. There is a formally smooth restriction map (see Remark \ref{rem-ChtR-product} for details)
\begin{equation}  \label{equation-R-nx-intro}
\Res^{nx} : (\Cht_{I \cup J \cup \{0\}, W \boxtimes V \boxtimes W_0})_{|\check{X}^{I \cup J} \times x} \to (L^+_{I \cup J}G \backslash \Gr_{I \cup J,W \boxtimes V})_{\check{X}^{I \cup J}} \times (\ChtR_{\{0\}, W_0, m, 0}^{nx})_x.
\end{equation} 
Let $\mcA \in \DD^b_c((\ChtR_{\{0\}, W_0, m, 0}^{nx})_x, \Lambda)$. Denote by $N^x = N - x$, we define 
$$\mcF_{I,W} = \mcL_{N^x} \otimes (\Res^{nx})^*(\mcS_{I \cup J, W \boxtimes V} \boxtimes \mcA).$$
If we need to specify the level structure, we will write $\mcF_{I,W,N^x}$. 
In the following, we will take nearby cycles for legs indexed by $I$, and do nothing for legs indexed by $J$.

\subsection{Reminder on nearby cycles}  \label{subsectionReminderNearbyCycles}

\subsubsection{}

Let $\mcY$ be a scheme (or stack) of finite presentation over a base $S$. For $t, s$ geometric points of $S$ and a specialization map $t \to s$ in $S$, we get a nearby cycle functor 
$$\Psi_{t \to s} : \DD_{\indcons}(\mcY_t, \Lambda) \to \DD_{\indcons}(\mcY_s, \Lambda)$$
defined as follows. First, consider the diagram 
\[\begin{tikzcd}
	{\mcY_t} & {\mcY_{S_{(s)}}} & {\mcY_s} \\
	t & {S_{(s)}} & s
	\arrow["j"{description}, from=1-1, to=1-2]
	\arrow["p"{description}, from=1-1, to=2-1]
	\arrow["p"{description}, from=1-2, to=2-2]
	\arrow["i"{description}, from=1-3, to=1-2]
	\arrow["p"{description}, from=1-3, to=2-3]
	\arrow[from=2-1, to=2-2]
	\arrow[from=2-3, to=2-2]
\end{tikzcd}\]
where $S_{(s)}$ is the strict henselization of $S$ at $s$ and $t \to S_{(s)}$ is our given specialization map. The nearby cycle functor is then defined as $$\Psi_{t \to s} = i^*j_*.$$
More generally given a sheaf $\mcF$ on $\mcY$ we will still denote by $\Psi_{t \to s}\mcF$ the sheaf obtained as the composition of the pullback along $\mcY_t \to \mcY$ and then $\Psi_{t \to s}$. 

\begin{rques}
\begin{enumerate}
\item In the case when $\mcY = S$, for all $\mcF \in \DD_c^b(\mcY, \Lambda)$ there is a canonical isomorphism $\Psi_{t \to s} \mcF = \mcF_t$ between the nearby cycle functor and the fiber at $t$. 
\item In the case when $S$ is of dimension one, $S_{(s)}$ is a trait, $\Psi_{t \to s}$ defined above is the classical nearby cycle functor for a trait (\cite[7, Exposé XIII]{SGA7}). 
\item In general, the above definition of $\Psi_{t \to s}$ computes the cohomology of Milnor fibers. We refer to Remark \ref{rqueOrgogozo} for the relation with the definition of nearby cycles over a general base (by Deligne, Laumon, Gabber, Orgogozo), which compute the cohomology of Milnor tubes. 
\item If we want to put some emphasis on the base for the calculation of the nearby cycles, we will denote by $\Psi^S_{t \to s}$ the functor $\Psi_{t \to s}$. In general, even if $t$ and $s$ both factor through some closed subscheme $Z \subset S$ the two functors $\Psi^Z_{t \to s}$ and $\Psi^S_{t \to s}$ may differ. This problem is resolved using the notion of $\Psi$-good sheaves, we refer to the appendix \ref{AppendixCyclesProches} for a discussion. 
\end{enumerate}
\end{rques}

\subsubsection{} 
In general, there is always a canonical base change map 
$$p_!\Psi_{t \to s} \to \Psi_{t \to s}p_!$$
coming from base change map $p_! i^* = i^* p_!$ and $p_!j_* \to j_*p_!$ (which is $p_!  j_* \xrightarrow{\mathrm{unit}} j_* j^* p_!  j_* \simeq j_* p_! j^* j_* \xrightarrow{\mathrm{counit}} j_* p_!$). 

When $p$ is a proper morphism, the above morphism of functor $p_!\Psi_{t \to s} \to \Psi_{t \to s}p_!$ is an isomorphism. In fact, in this case, $p_! = p_*$, we have $p_*  j_* \simeq j_* p_*$.

\subsubsection{}
In the specific case where $S = X$, we will denote by $\Psi$ the nearby cycle functor along the specialization map $\mathfrak{sp}: \overline{\eta} \to x$ that we have already chosen. More generally if $S = X^I$, we will denote by $\Psi_i$ the nearby cycle functor defined by the map $\mcY \to S \xrightarrow{\pr_i} X$ by the $i$th projection and the specialization $\mathfrak{sp}$.

\subsection{Main statements}

We fix $(J,V)$ and $\mcA$ as above. 

\begin{thm}\label{thmMain2}
For all totally ordered $I = \{1, \dots, n\}$, and $W \in \Rep_{\Lambda}{^L}G^I$, the canonical map 
\begin{equation}\label{equationMapTheorem2}
\mathfrak{p}_!\Psi_1\dots\Psi_n(\mcF_{I,W}) \rightarrow \Psi_1\dots\Psi_n \mathfrak{p}_!(\mcF_{I,W})
\end{equation}
is an isomorphism. 
\end{thm}

We denote by $\eta_I$ the generic point of $X^I$ and by $\overline{\eta}_I$ a geometric point over $\eta_I$. We fix some specialization map $\overline{\eta}_I \to x$.

\begin{thm}\label{thmMain1}
\begin{enumerate}
\item Assume that $\Lambda$ is finite. We denote by $\Delta : X \to X^I$ the diagonal map. For all $(I,W)$, there exists a modification $\widetilde{S} \to X^I$,  a point $s \in \widetilde{S} \times_{X^I} \Delta(x)$ and a specialization map $\overline{\eta}_I \to s$ such that the canonical map
\begin{equation}\label{equationMapTheorem1}
\mathfrak{p}_!\Psi_{\overline{\eta}_I \to s}(\mcF_{I,W}) \rightarrow \Psi_{\overline{\eta}_I \to s} \mathfrak{p}_! (\mcF_{I,W})
\end{equation}
is an isomorphism. 

\item Moreover the point $s$ can be chosen such that there exists a canonical isomorphism of sheaves on $(\Cht_{I \cup J \cup \{0\} })_{| x^I \times \check{X}^{J} \times x}$
\begin{equation*}
\Psi_{\overline{\eta}_I \to s}(\mcF_{I,W}) = \Psi_1\dots\Psi_n(\mcF_{I,W}). 
\end{equation*}
Besides we also have a canonical isomorphism of sheaves on $\check{X}^J$
\begin{equation*}
\Psi_{\overline{\eta}_I \to s} \mathfrak{p}_!(\mcF_{I,W}) = \Psi_1\dots\Psi_n \mathfrak{p}_!(\mcF_{I,W}).
\end{equation*}
The following diagram is commutative:
\[\begin{tikzcd}
	{\mathfrak{p}_!\Psi_{\overline{\eta}_I \to s}(\mcF_{I,W})} & {\Psi_{\overline{\eta}_I \to s} \mathfrak{p}_!(\mcF_{I,W})} \\
	{p_!\Psi_1\dots\Psi_n(\mcF_{I,W})} & {\Psi_1\dots\Psi_n \mathfrak{p}_!(\mcF_{I,W}).}
	\arrow["(\ref{equationMapTheorem1})", from=1-1, to=1-2]
	\arrow["\simeq", from=1-1, to=2-1]
	\arrow["\simeq", from=1-2, to=2-2]
	\arrow["(\ref{equationMapTheorem2})", from=2-1, to=2-2]
\end{tikzcd}\]
\end{enumerate}
\end{thm}

In section \ref{sectionOneLeg}, we illustrate the argument in the case where $I = \{1\}, N^x = \emptyset$ and $\mcA = \mcL_{nx}$ which is the simplest possible case. In section \ref{sectionFusion}, we prove some technical results needed in section \ref{sectionGeneralCase}, where we prove the main theorems.

\subsection{Motivation and relation with litteratures}

\subsubsection{}  \label{subsection-lien-avec-Salmon}
Such type of statements appeared in \cite{Salmon}, \cite{Salmon2} and in the thesis of the first author \cite{EteveThese}. In \cite{Salmon2} the author shows theorem \ref{thmMain2} under an extra assumption which he calls $\Psi$-factorizability. This assumption can be thought of as an assumption on the sheaf $\mcA$, here we completely remove this assumption at the cost of working at the level of stacks of shtukas (and restricted shtukas) contrary to Salmon who works at the level of affine Grassmannians and affine flag varieties. Our argument is a refinement of the previous ones and involves carefully constructing fusion isomorphisms for the nearby cycles that appear. In \cite{Salmon2} and \cite{EteveThese}, the authors use theorem \ref{thmMain2} to compute the monodromy of the cohomology of shtukas when the level is reduced a place, we expect that our present result will be of use to control the monodromy for deeper levels.

\subsubsection{}
To simplify the notation let $J$ be empty.
Theorem \ref{thmMain1} is also true if we replace $x \in N$ by some $u \in \check{X}^I$. This implies that $\mathfrak{p}_! \mathcal{F}_{I, W}$ is ind-lisse over $\check{X}^I$. In fact, for any geometric point $u \in \check{X}^I$ and any specialization map $\mathfrak{sp}_u: \overline{\eta}_I \rightarrow u$, we have a canonical commutative diagram:
$$
\xymatrix{
\mathfrak{p}_!\Psi_{\overline{\eta}_I \to u}(\mcF_{I,W})  \ar[rr]^{(\ref{equationMapTheorem1})}  
& & \mathfrak{p}_! (\mcF_{I,W}) |_{ \overline{\eta}_I }  \\
& \mathfrak{p}_! (\mcF_{I,W}) |_{u}   \ar[lu]^{\mathrm{adj}}  \ar[ru]_{\mathfrak{sp}_u^*}
}
$$
where $\mathrm{adj}$ is induced by $i^* \mathcal{F}_{I, W} \rightarrow i^* j_*j^* \mathcal{F}_{I, W}$. By the geometric Satake equivalence, $\mathcal{S}_{I, W}$ is ULA. We deduce that $\mathcal{F}_{I, W}$ is also ULA. So morphism $\mathrm{adj}$ is an isomorphism. The statements ``$(\ref{equationMapTheorem1})$ is an isomorphism" and ``$\mathfrak{sp}_u^*$ is an isomorphism" are equivalent.

In particular, when $\mathcal{A} = \mathcal{L}_{nx}$, the complex $\mathfrak{p}_! \mathcal{F}_{I, W}$ is the cohomology of shtukas $\mathcal{H}_{G, N, I, W}$ considered in \cite{XueSmoothness}. We refind the result that $\mathcal{H}_{G, N, I, W}$ is ind-lisse over $\check{X}^I$. (The proofs are the same: if we replace everywhere $\mathfrak{p}_!\Psi_{\overline{\eta}_I \to u}(\mcF_{I,W})$ by $\mathfrak{p}_! (\mcF_{I,W}) |_{u}$ in the ``Zorro lemma" argument, the diagrams in section \ref{sectionGeneralCase} coincide with the diagrams in \cite{XueSmoothness}.)

%\subsubsection{}
%When $I$ is a singleton, Theorem \ref{thmMain2} just says that the classical nearby cycle functor commutes with proper direct image functor (for the sheaf $\mcF$). In other words, the canonical morphism
%\begin{equation}    \label{equation-H-c-Cht-x-to-H-c-Cht-eta} 
%\mathrm{can}: H_c(\Cht | _{x}, \Psi \mathcal{F}) \rightarrow H_c(\Cht |_{\overline{\eta}}, \mathcal{F})
%\end{equation}
%is an isomorphism.
%Note that this morphism is equivarient under the Galois action, for action of $\mathrm{Gal}(\overline{F} / F)$ on the cohomology group on the RHS and the action of $\mathrm{Gal}(\overline{F}_x / F_x)$ on the sheaf $\Psi \mathcal{F}$ on the LHS.

\subsubsection{}  \label{subsection-local-global-compatibility}
One motivation of this paper is that we expect the main result theorem \ref{thmMain2}, which is of technical nature, to be the main input in showing a strong form of local-global compatibility. We refer to \cite[Section 4.5]{HowToInventShtukas} for a discussion. The methods used here will also be used in \cite{LocalSpectralAction} to construct a spectral action comparable to the one of \cite{FarguesScholze}.

\subsubsection{}
The map (\ref{equation-R-nx-intro}) factors through 
$\Cht_{I \cup \{0\}, W \boxtimes W_0} \to \ChtR_{I \cup \{0\}, W \boxtimes W_0}^{nx}$, which is also a smooth morphism. Since smooth pullback commutes with nearby cycles, the nearby cycles we considered in this paper (over stacks of global shtukas) are the same as the ones in \cite{GenestierLafforgue} (over stacks of restricted shtukas), where the authors consider the case $\mathcal{A} = \mathcal{L}^{nx}$.

\subsection{Reduction to the modular case.}

In the rest of this paper, we will need to apply some strong theorems about higher dimensional nearby cycles which \emph{a priori} only hold for torsion coefficients. We therefore want to indicate how to reduce the proof of theorem \ref{thmMain2} to the case where $\Lambda$ is a finite field. 

Assume that $\Lambda = E$, since the map (\ref{equationMapTheorem2}) is functorial we can proceed by devissage and assume that $\mcA = \mcA' \otimes_{\mcO_E} E$ for some sheaf $\mcA' \in \DD^b_c((\ChtR_{\{0\}, W_0, m, 0}^{nx})_x, \mcO_E)$. Similarly, we can assume that $W \boxtimes V = (W' \boxtimes V') \otimes_{\mcO_E} E$ for a representation $(W' \boxtimes V')$ of ${^L}G^{I \cup J}$ on a finite free $\mcO_E$ module. Since both functors $\Psi_i$ and $\mathfrak{p}_!$ commute with $\otimes_{\mcO_E} E$ it is enough to show that the map (\ref{equationMapTheorem2}) is an isomorphism when $\Lambda = \mcO_E$. Let $C_{\mcO_E}$ be the cone of the map (\ref{equationMapTheorem2}) since both functor $\Psi_i$ and $\mathfrak{p}_!$ also commute with reduction mod $\ell$, we have that 
$$C_{\mcO_E} \otimes_{\mcO_E} k_E = \mathrm{cone}\big( \mathfrak{p}_!\Psi_{\overline{\eta}_I \to s}(\mcF_{I,W} \otimes_{\mcO_E} k_E) \rightarrow \Psi_{\overline{\eta}_I \to s} \mathfrak{p}_! (\mcF_{I,W} \otimes_{\mcO_E} k_E) \big).$$
Assuming that theorem \ref{thmMain2} holds if $\Lambda = k_E$, we deduce that $ C_{\mcO_E} \otimes_{\mcO_E} k_E = 0$. Using the triangle 
$$C_{\mcO_E} \xrightarrow{\varpi} C_{\mcO_E} \to C_{\mcO_E} \otimes_{\mcO_E} k_E$$
where the first map is the multiplication by a uniformizer $\varpi$ of $\mcO_E$, we get short exact sequences 
$$0 \to H^i(C_{\mcO_E})/\varpi \to H^i(C_{\mcO_E} \otimes_{\mcO_E} k_E) \to \mathrm{Tor}^1_{\mcO_E}(H^{i+1}(C_{\mcO_E}), k_E) \to 0,$$ 
as the middle term vanishes so do the first and last ones. Let $u$ be an unramified place for $G$ and not in $N_0$, then by a mild generalization of 
\cite[Theorem 0.0.3]{XueIntegral}, the cohomology groups of $\mathfrak{p}_! \mcF_{I,W}$ and $\mathfrak{p}_!\Psi_{\overline{\eta}_I \to s}\mcF_{I,W}$ are modules of finite type over the local Hecke algebra at $u$. Consequently, the cohomology groups of $C_{\mcO_E}$ are also of finite type over this Hecke algebra, by Nakayama's lemma we deduce that $C_{\mcO_E}^i = 0$ and thus that $C_{\mcO_E} = 0$. 

\begin{rque}
Note that $\mathfrak{p}_! \mcF_{I,W}$ and $\mathfrak{p}_!\Psi_{\overline{\eta}_I \to s}\mcF_{I,W}$ are complex of ind-constructible $\mcO_E$-sheaves. The above discussion shows that they do not contain a copy $E$ as ind-constructible $\mcO_E$-sheaves which would contradict the application of Nakayama's lemma.
\end{rque}

\begin{center}
\textbf{From now on, we assume that $\Lambda$ is a torsion ring killed by a power of $\ell$.}
\end{center}

\section{Case of one leg}\label{sectionOneLeg}

The goal of this section is to give a proof of Theorem \ref{main-thm-one-leg}, which illustrates the proof of general case. 
In \ref{subsection-relation-with-mainthms}, we explain why Theorem \ref{main-thm-one-leg} is a special case of Theorem \ref{thmMain2}. 
In the rest of this section, we prove Theorem \ref{main-thm-one-leg}. We will first construct a morphism in the inverse direction of (\ref{equation-p-Psi-F-to-Psi-p-F-one-leg}) in \ref{subsectionInversemap}-\ref{subsection-ourcase-oneleg}, then show that it is indeed an inverse of (\ref{equation-p-Psi-F-to-Psi-p-F-one-leg}) in \ref{subsectionZorro-oneleg}.

\subsection{Relation of Theorem \ref{main-thm-one-leg} and Theorem \ref{thmMain2}}  \label{subsection-relation-with-mainthms}

When $I$ is a singleton, Theorem \ref{thmMain2} and Theorem \ref{thmMain1} are the same. In the following we explain why Theorem \ref{main-thm-one-leg} is a special case of Theorem \ref{thmMain2}.

\subsubsection{}
In Theorem \ref{thmMain2}, let $I=\{1\}$ and $J$ be the empty set. 
Let $N^x$ be the empty set. In this case, $N=nx$. We also suppose that $G$ is split and constant on $X$ (in particular, $R$ is empty). Also suppose that $\deg(x)=1$. We have $\check{X}=X-x$.

The smooth morphism (\ref{equation-R-nx-intro}) is
\begin{equation}   \label{equation-R-nx-one-leg}
\mathcal{R}^{nx}: (\Cht_{\{1, 0\}, W \boxtimes W_0}) |_{(X -x) \times x} \rightarrow [L^+_{\{1\}}G \backslash \Gr_{\{1\}, W}] |_{X -x} \times \ChtR^{nx}_{\{0\}, W_0, m} |_x
\end{equation}
For any $\mathcal{A} \in D_c^b(\Cht\mathcal{R}^{nx}_{\{0\}, W_0, m} |_{x}, \Lambda)$, the sheaf $\mcF_{\{1\}, W}$ is the pullback
$(\mathcal{R}^{nx})^*(\mathcal{S}_{\{1\}, W} \boxtimes \mathcal{A})$  over $(\Cht_{\{1, 0\}, W \boxtimes W_0}) |_{(X -x) \times x}$.

\subsubsection{}
We denote by $\mathbf{1}$ the trivial representation of ${^L}G$.
Recall that we have
$$\Cht_{\{1, 0\}, W \boxtimes \mathbf{1}} = \Cht_{\{1\}, W} \times X.$$

For $W_0 = \mathbf{1}$, we have $\ChtR^{nx}_{\{0\}, \mathbf{1}} = \bullet / G(\mathcal{O}_{nx})$. 
Recall that $N=nx$. We have the following Cartesian diagram:
$$
\xymatrix{
(\Cht_{\{1\}, W, N}) |_{X -x} \times \bullet  \ar[rr]  \ar[d]^{\pi}
& & [L^+_{\{1\}}G \backslash \Gr_{\{1\}, W}]|_{X-x} \times \bullet \ar[d] \\
(\Cht_{\{1\}, W}) |_{X -x} \times \bullet  \ar@{=}[r]  
& (\Cht_{\{1, 0\}, W \boxtimes \mathbf{1}}) |_{(X -x) \times x}   \ar[r]^{\mathcal{R}^{nx} \quad \quad }
& [L^+_{\{1\}}G \backslash \Gr_{\{1\}, W}]|_{X-x} \times \bullet / G(\mathcal{O}_{nx})
}
$$
where the objects in the upper line are $G(\mathcal{O}_{nx})$-torsors over the objects in the lower line.
Let $\mathcal{A}$ over $\bullet / G(\mathcal{O}_{nx})$ be the direct image of $\Lambda$ via $\bullet \rightarrow \bullet / G(\mathcal{O}_{nx})$. 
We have $$(\mathcal{R}^{nx})^*(\mathcal{S}_{\{1\}, W} \boxtimes \mathcal{A}) = \mathcal{S}_{\{1\}, W} \otimes \mathcal{L}^{nx}$$ over $(\Cht_{\{1\}, W}) |_{X -x}$, where $\mathcal{L}^{nx} = \pi_! \Lambda$ is viewed as a sheaf over $(\Cht_{\{1\}, W}) |_{X -x}$.

\subsubsection{}
The nearby cycle functor in Theorem \ref{thmMain2} are taken for the following commutative diagram 
$$
\xymatrix{
(\Cht_{\{1, 0\}, W \boxtimes W_0})|_{\overline{\eta} \times x}  \ar[r]^{j}  \ar[d]^{\mathfrak{p}}
& \Cht_{\{1, 0\}, W \boxtimes W_0}  \ar[d]^{\mathfrak{p}}
& (\Cht_{\{1, 0\}, W \boxtimes W_0})|_{x \times x }   \ar[l]_{i}  \ar[d]^{\mathfrak{p}} \\
\overline{\eta} \times x \ar[r]^{j}
& X_{(x)} \times x
& x \times x  \ar[l]_{i}
}
$$
All the products are taken over $\mathrm{Spec}(\Fqb)$. So we can identify $\overline{\eta} \times x$ with $\overline{\eta}$ and identify $x \times x$ with $x$. When $W_0 = \mathbf{1}$, this diagram coincides with the diagram in \ref{subsection-nearby-cycle-one-leg}. 

\quad

Thus Theorem \ref{main-thm-one-leg} is a special case of Theorem \ref{thmMain2}, for $I=\{1\}$, $J$ be the empty set, $W_0 = \mathbf{1}$ the trivial representation and $\mathcal{A} = \mathcal{L}^{nx}$.

\subsection{Construction of the inverse map}  \label{subsectionInversemap}

\subsubsection{}
To construct the inverse map, we need stacks of shtukas with several legs. For any $(I', W')$, we denote by $\mathcal{F}_{I', W'} := \mathcal{S}_{I', W'} \otimes \mathcal{L}^{nx}.$
Let $\gamma$ be the following composition of morphisms:
\begin{equation}   \label{equation-inverse-of-can}
\begin{aligned}
\Psi_1 \mathfrak{p}_! \mathcal{F}_{\{1\}, W}  \otimes \Lambda  & \simeq \Psi_1 \mathfrak{p}_! \Psi_2 \mathcal{F}_{\{1, 2\}, W \boxtimes \mathbf{1}} \\
& \xrightarrow{\delta} \Psi_1 \mathfrak{p}_! \Psi_2 \mathcal{F}_{\{1, 2\}, W \boxtimes (W^* \otimes W)}  \\
& \simeq \Psi_1 \mathfrak{p}_!  \Psi_2 (\Delta^{\{2, 3\}})^* \mathcal{F}_{\{1, 2, 3\}, W \boxtimes W^* \boxtimes W}  \\
& \xrightarrow{\simeq \; \alpha} \Psi_1 \mathfrak{p}_!  \Psi_2 \Psi_3 \mathcal{F}_{\{1, 2, 3\}, W \boxtimes W^* \boxtimes W}  \\
& \xrightarrow{\mathrm{can}} \Psi_1 \Psi_2 \mathfrak{p}_!   \Psi_3 \mathcal{F}_{\{1, 2, 3\}, W \boxtimes W^* \boxtimes W}  \\
& \xrightarrow{\simeq \; \beta^{-1}} \Psi_2 (\Delta^{\{1, 2\}})^* \mathfrak{p}_!  \Psi_3 \mathcal{F}_{\{1, 2, 3\}, W \boxtimes W^* \boxtimes W}  \\
& \simeq \Psi_2 \mathfrak{p}_! \Psi_3 \mathcal{F}_{\{2, 3\}, (W \otimes W^*) \boxtimes W}  \\
& \xrightarrow{\mathrm{ev}} \Psi_2 \mathfrak{p}_! \Psi_3 \mathcal{F}_{\{2, 3\}, \mathbf{1} \boxtimes W}  \\
& \simeq \Lambda \otimes \mathfrak{p}_! \Psi_3 \mathcal{F}_{\{3\}, W}
\end{aligned}
\end{equation}

\quad

The maps are the following ones:

(1) The first isomorphism is because 
$$\Cht_{\{1, 2\}, W \boxtimes \mathbf{1}} = \Cht_{\{1\}, W} \times X$$
and the sheaf $\mathcal{F}_{\{1, 2\}, W \boxtimes \mathbf{1}}$ is isomorphic to $\mathcal{F}_{\{1\}, W} \boxtimes \Lambda$.
Thus $$\Psi_1 \mathfrak{p}_! \Psi_2 \mathcal{F}_{\{1, 2\}, W \boxtimes \mathbf{1}} \simeq \Psi_1 \mathfrak{p}_! \Psi_2 (\mathcal{F}_{\{1\}, W} \boxtimes \Lambda) = \Psi_1 \mathfrak{p}_! (\mathcal{F}_{\{1\}, W} \otimes \Lambda) = \Psi_1 \mathfrak{p}_! \mathcal{F}_{\{1\}, W} \otimes \Lambda$$

(2) The second morphism follows from the functoriality of Satake sheaves associated to the canonical morphism $\delta: \mathbf{1} \rightarrow W^* \otimes W$.

(3) The third isomorphism follows from the fusion of Satake sheaves associated to $\{1, 2, 3\} \rightarrow \{1, 2\}$ sending $1$ to $1$ and $\{ 2, 3\}$ to $2$, where $\Delta^{\{2, 3\}}: X^2 \rightarrow X^3, (x_1, x_2) \mapsto (x_1, x_2, x_2)$ is the partial diagonal inclusion.

The composition of the first three morphisms is $\Psi_1 \mathfrak{p}_!  \Psi_2$ applied to the following creation morphism (creats legs $2$ and $3$) of \protect{\cite[Section 5]{Lafforgue}}:
$$\mathcal{C}^{\sharp}_{\{2, 3\}}:  \mathcal{F}_{\{1\}, W}  \otimes \Lambda \rightarrow  (\Delta^{\{2, 3\}})^* \mathcal{F}_{\{1, 2, 3\}, W \boxtimes W^* \boxtimes W}$$

\quad

(4) The fourth isomorphism $\alpha$ is the key morphism and is the most difficult one to construct. It will be given in Construction \ref{constr-morphism-alpha}.

(5) The fifth morphism is the canonical morphism $\mathfrak{p}_! \Psi_2 \rightarrow \Psi_2 \mathfrak{p}_!$.

(6) The sixth isomorphism $\beta^{-1}$ will be constructed in Construction \ref{constr-morphism-beta}. The construction uses Proposition \ref{prop-p-F-constant}, whose proof uses Drinfeld's lamma.

\quad

(7) The seventh isomorphism follows from the fusion of Satake sheaves associated to $\{1, 2, 3\} \rightarrow \{2, 3\}$ sending $\{1, 2 \}$ to $2$ and $3$ to $3$, where $\Delta^{\{1, 2\}}: X^2 \rightarrow X^3, (x_1, x_2) \mapsto (x_1, x_1, x_2)$ is the partial diagonal inclusion.

(8) The eighth morphism follows from the functoriality of Satake sheaves associated to the evaluation map $\mathrm{ev}: W \otimes W^* \rightarrow \mathbf{1}$.

(9) The last isomorphism is because 
$$\Cht_{\{1, 3\}, \mathbf{1} \boxtimes W} =X \times  \Cht_{\{3\}, W}$$
and the sheaf $\mathcal{F}_{\{1, 3\}, \mathbf{1} \boxtimes W}$ is isomorphic to $\Lambda \boxtimes \mathcal{F}_{\{3\}, W}$.
Thus $$\Psi_1 \mathfrak{p}_! \Psi_3 \mathcal{F}_{\{1, 3\}, \mathbf{1} \boxtimes W } \simeq \Psi_1 \mathfrak{p}_! \Psi_3 (\Lambda \boxtimes \mathcal{F}_{\{3\}, W}) =\Lambda \otimes \mathfrak{p}_! \Psi_3 \mathcal{F}_{\{3\}, W} . $$

The composition of the last three morphisms is $\Psi_2  \mathfrak{p}_!  \Psi_3$ applied to the following annihilation morphism (annihilates legs $1$ and $2$) of \protect{\cite[Section 5]{Lafforgue}}:
$$\mathcal{C}^{\flat}_{\{1, 2\}}: (\Delta^{\{1, 2\}})^*  \mathcal{F}_{\{1, 2, 3\}, W \boxtimes W^* \boxtimes W}  \rightarrow \Lambda \otimes  \mathcal{F}_{\{3\}, W}$$

\subsubsection{}
To construct morphism $\alpha$, we need some commutativity of nearby cycles with the partial diagonal restriction. 
Let's first consider the simplest case, when there is only two legs, in \ref{subsectionToycase}. Then we treat our case of three legs in \ref{subsection-ourcase-oneleg}.

\begin{rque}
The principal difficulty in the construction of (\ref{equation-inverse-of-can}) is to construct some kind of isomorphism of type $\Psi \Delta^* \mcF = \Psi_1 \Psi_2 \mcF$. This is true if $\mcF$ is of the form $\mathcal{F}_1 \boxtimes \mathcal{F}_2$ over a scheme or stack of the form $\mathcal{Y}_1 \times \mathcal{Y}_2$ over $X^2$ (because of Künneth formula).

In our case, when $\mathcal{A} = \Lambda$ is the constant sheaf, we are in the situation of (\ref{equation-epsilon}), we can reduce the calculation to the product of affine grassmannians and apply the Künneth formula (because the product structure is over $X^2$).

For general $\mathcal{A}\in \DD^b_c((\ChtR_{\{0\}, W_0}^{nx})_x, \Lambda)$, we are in the situation of (\ref{equation-R}), the product structure is only over $(X-x)^2 \times x$, and we do not have a product structure over $X^2 \times x$. So we can not apply the Künneth formula. To solve this problem, we use the twisted product structure of stacks of restricted shtukas, and the partial Frobenius morphism, to get a weaker statement: $\mathfrak{p}_! \Psi \Delta^* \mcF = \mathfrak{p}_! \Psi_1 \Psi_2 \mcF$. This is explained in the section \ref{subsectionToycase} below.
\end{rque}

\subsection{Toy case: nearby cycle commutes with the diagonal restriction}
\label{subsectionToycase}

\subsubsection{Nearby cycles on shtukas}

Let $I=\{1, 2\}$ and $\mathfrak{p}: \Cht_{\{1, 2\}} \rightarrow X^2$. Let $\Delta: X \rightarrow X^2$ be the diagonal inclusion. Let $$\mathrm{Frob}_{\{1\}}: X^2 \rightarrow X^2, (x_1, x_2) \mapsto (\mathrm{Frob}(x_1), x_2)$$ be the partial Frobenius morphism, where $\mathrm{Frob} : X \to X$ is the $q$-power $\Fqb$-linear Frobenius. For any $d \in \mathbb{Z}_{\geq 0}$, let $\Delta_d = \mathrm{Frob}_{\{1\}}^d \Delta$. 
For simplicity suppose that $\mathrm{deg}(x)=1$ (if not, use $\mathrm{Frob}_{\{1\}}^{\mathrm{deg(x)}}$ instead of $\mathrm{Frob}_{\{1\}}$). We have $\Delta_d(x)=\Delta(x)$. Let $\mathcal{F} = \mathcal{S}_{\{1, 2\}, W_1 \boxtimes W_2} \otimes \mathcal{L}^{nx}$. 

\begin{lem}   \label{lem-two-legs-Psi-delta-d-eta-equal-Psi-1-Psi-2}
For $d \gg 0$, we have a canonical morphism 
\begin{equation}   \label{equation-Psi-Delta-d-to-Psi-1-Psi-2}
\alpha: \Psi \Delta_d^* \mathcal{F} \rightarrow \Psi_1 \Psi_2 \mathcal{F}
\end{equation}
which is an isomorphism.
\end{lem}

Even if the statement of this lemma involves only usual nearby cycles over a basis of dimension one, the proof needs nearby cycles over a general basis. 

\begin{proof}
By Orgogozo's theorem (Theorem \ref{thmOrgogozo}, Remark \ref{rqueOrgogozo}), there exists a modification $\widetilde{S} \rightarrow X^2$ such that (a) over $\widetilde{S}$ the sheaf $\mathcal{F}$ is $\Psi$-good (i.e. the formation $\Psi_{*} \mathcal{F}$ commutes with all base changes), 
(b) %for any specialization map $a \rightarrow b$ in $\widetilde{S}$, the restriction of $\Psi_{*} \mathcal{F}$ to $\Cht_{\{1,2\}} |_b$ is isomorphic to $\Psi_{a \rightarrow b} \mathcal{F}$ 
the definition $\Psi_* \mcF$ coincides with the naive definition in \ref{subsectionReminderNearbyCycles} (i.e. cohomology of Milnor tubes equals to cohomology of Milnor fibers). 

For the flag $x \subset X \times x \subset X^2$, we have a sequence of strict transforms (see \ref{subsection-stricttransform} for the definition):
\begin{equation}  \label{equation-flag-hyperplan}
\xymatrix{
s  \ar@{^{(}->}[r]  \ar[d]
& \widetilde{S}_1 \ar@{^{(}->}[r]  \ar[d]
& \widetilde{S} \ar[d] \\
x \ar@{^{(}->}[r]
& X \times x \ar@{^{(}->}[r] 
& X^2
}
\end{equation}
For any $d \in \mathbb{Z}_{\geq 0}$, we have a flag $x \subset \Delta_d(X) \subset X^2$ and a sequence of strict transforms:
\begin{equation}   \label{equation-flag-diagonal}
\xymatrix{
s_d  \ar@{^{(}->}[r]  \ar[d]
& \widetilde{S}_{1, d} \ar@{^{(}->}[r]   \ar[d]
& \widetilde{S} \ar[d] \\
x \ar@{^{(}->}[r]
& \Delta_d(X) \ar@{^{(}->}[r]
& X^2
}
\end{equation}
Let $d \in \mathbb{Z}_{\geq 0}$ be large enough such that  
$s_d = s$. 
We view $\overline{\eta_I}$ as a geometric generic point of $\widetilde{S}$.
Consider (\ref{equation-flag-hyperplan}), let $\overline{\eta}_1$ be a geometric point of $\widetilde{S}_1$ over $\overline{\eta} \times x$. 
Consider (\ref{equation-flag-diagonal}), we view $\Delta_d(\overline{\eta})$ as a geometric point of $\widetilde{S}_{1, d}$.
Choose specialization maps in $\widetilde{S}$ such that the following diagram is commutative:
$$
\xymatrix{
\overline{\eta_I}  \ar[r]  \ar[d]
& \Delta_d(\overline{\eta}) \ar[d] \\
\overline{\eta}_1  \ar[r]
& s=s_d
}
$$

We construct canonical isomorphisms:
\begin{equation*}
\begin{aligned}
\Psi \Delta_d^* \mathcal{F}  & \simeq \Psi_{\Delta_{d}(\overline{\eta}) \rightarrow s_d} \mathcal{F}  \\
& \simeq \Psi_{\overline{\eta}_I \rightarrow s_d} \mathcal{F} \\
& \simeq \Psi_{\overline{\eta}_I \rightarrow s} \mathcal{F} \\
& \simeq \Psi_{\overline{\eta}_1 \rightarrow s} \Psi_{\overline{\eta}_I \rightarrow \overline{\eta}_1} \mathcal{F}  \\
& \simeq \Psi_1 \Psi_2 \mathcal{F}
\end{aligned}
\end{equation*}
The maps are the following:

\begin{enumerate}
\item the first isomorphism comes from the base change map (here we use Orgogozo's theorem (a) and (b), which is also Remark \ref{rquebasechange})
$$\Psi_{\Delta_{d}(\overline{\eta}) \rightarrow s_d} \mathcal{F}  \xrightarrow{\sim}  \Psi \Delta_d^* \mathcal{F}  $$

\item the second isomorphism comes from Gabber's theorem (Theorem \ref{thmGabber}) for the sequence of specialization maps $\overline{\eta}_I \rightarrow  \Delta_{d}(\overline{\eta}) \rightarrow s_d$:
$$\Psi_{\overline{\eta}_I \rightarrow s_d} \mathcal{F} \xrightarrow{\sim}  \Psi_{\Delta_{d}(\overline{\eta}) \rightarrow s_d} \Psi_{\overline{\eta}_I \rightarrow \Delta_{d}(\overline{\eta})} \mathcal{F} $$
Moreover both $\overline{\eta}_I$ and $\Delta_{d}(\overline{\eta})$ lie in $(X-x)^2$, where $\mathcal{F}$ is ULA. So $\Psi_{\overline{\eta}_I \rightarrow \Delta_{d}(\overline{\eta})} \mathcal{F} \simeq \mathcal{F}$. 

\item the third isomorphism is because $s_d=s$ and the commutativity of specialization maps.

\item the fourth isomorphism again comes from Gabber's theorem (Theorem \ref{thmGabber}) for the sequence of specialization maps $\overline{\eta}_I \rightarrow \overline{\eta}_1 \rightarrow s$:
$$\Psi_{\overline{\eta}_I \rightarrow s} \mathcal{F} \xrightarrow{\sim}   \Psi_{\overline{\eta}_1\rightarrow s} \Psi_{\overline{\eta}_I \rightarrow \overline{\eta}_1} \mathcal{F} $$

\item the last isomorphism comes from the Künneth formula (Theorem \ref{thmKunnethCyclesProches}) and the structure of twisted product of stacks of restricted shtukas. 

Let us detail this last point: take the notations in \ref{subsection-ChtR-intermediate}. Let $I_1 = \{1\}$ and $I_2 = \{2\}$. We want to construct the desired morphism over $\Cht_{\{1, 2\}}$. Since the morphism $\mathcal{R}^{nx}$ is smooth, $(\mathcal{R}^{nx})^*$ commutes with $\Psi$, it is enough to construct the desired morphism over $\ChtR_{\{1, 2\}}^{nx}$. Moreover, since the morphism $\pi$ is proper, $\pi_!$ commutes with $\Psi$, it is enough to construct the desired morphism over $\ChtR^{nx, (1, 2)}_{\{1, 2\}}$. 
Over $(X - x) \times X$, the stack of restricted shtukas $\ChtR^{nx, (1, 2)}_{\{1, 2\}}$ has the structure of twisted product (Lemma \ref{lemTwistedProductStructure}),
our sheaf $\mathcal{S}_{\{1, 2\}, W_1 \boxtimes W_2}^{(1, 2)} \otimes \mathcal{L}^{nx}$ is the inverse image of the twisted product $ \mathcal{S}_{\{1\}, W_1} \tilde{\boxtimes} (\mathcal{S}_{\{2\}, W_2} \otimes \mathcal{L}^{nx})$. 

Apply the Künneth formula (Theorem \ref{thmKunnethCyclesProches}) to 
$$[G_{n_1y_1} \backslash \Gr_{\{1\}, W_1}] |_{\overline{\eta}}  \; \widetilde{\times} \; \ChtR_{\{2\}, W_2, n_2}^{nx}  \xrightarrow{(\mathrm{pr}_1, \mathrm{pr}_2)} \overline{\eta} \times X $$
Let $\mathcal{A}_1= \mathcal{S}_{\{1\}, W_1}$ and $\mathcal{A}_2 = \mathcal{S}_{\{2\}, W_2} \otimes \mathcal{L}^{nx}$. 
The Künneth formula says that the canonical morphism
\begin{equation}  \label{equation-Kunneth-eta-12-to-eta-1}
\Psi_{\overline{\eta} \rightarrow \overline{\eta}} \mathcal{A}_1 \boxtimes \Psi_{\overline{\eta} \rightarrow x} \mathcal{A}_2 \rightarrow  \Psi_{\overline{\eta}_I \rightarrow \overline{\eta}_1} (\mathcal{A}_1 \boxtimes \mathcal{A}_2)
\end{equation}
is an isomorphism. Note that the LHS
$= \mathcal{A}_1 \boxtimes \Psi_{\overline{\eta} \rightarrow x} \mathcal{A}_2 = \Psi_2 (\mathcal{A}_1 \boxtimes \mathcal{A}_2)$
where $\Psi_2 $ is the naive nearby cycle for the projection $\mathrm{pr}_2$. So (\ref{equation-Kunneth-eta-12-to-eta-1}) is
\begin{equation*}
\Psi_2 \mathcal{F} \xrightarrow{\sim} \Psi_{\overline{\eta}_I \rightarrow \overline{\eta}_1} \mathcal{F}
\end{equation*}
Applying $\Psi_{\overline{\eta}_1 \rightarrow s}$ to two sides,
we deduce
\begin{equation*}  
\Psi_1  \Psi_2 \mathcal{F} \xrightarrow{\sim}  \Psi_{\overline{\eta}_1 \rightarrow s}  \Psi_{\overline{\eta}_I \rightarrow \overline{\eta}_1 } \mathcal{F}
\end{equation*}
\end{enumerate}

\end{proof}

\begin{rque}
We can not use Künneth formula for
$$\Psi_{\overline{\eta} \rightarrow x} \mathcal{A}_1 \boxtimes \Psi_{\overline{\eta} \rightarrow x} \mathcal{A}_2 \rightarrow  \Psi_{\overline{\eta}_{I} \rightarrow s} (\mathcal{A}_1 \boxtimes \mathcal{A}_2)$$
because $\ChtR_{\{1, 2\}}^{nx, (1,2)}$ does not have a structure of twisted product over $X_{(x)} \times X_{(x)}$. It only has a structure of twisted product over $\overline{\eta} \times X_{(x)}$. That's why we have to combine Step (4) and Step (5) in the above proof.
\end{rque}

\begin{rque}
For general case of Step (2), see \ref{subsubsectionPartialDiagonal}. For general case of Step (4) and Step (5), see \ref{subsubsectionIteration}. 
\end{rque}

\begin{lem}   \label{lem-p-Frob-d-Delta-equal-p-Delta}
We take the same notations as in Lemma \ref{lem-two-legs-Psi-delta-d-eta-equal-Psi-1-Psi-2}. For any $d \in \mathbb{Z}_{\geq 0}$, we have $$\mathfrak{p}_! \Psi \Delta_d^*  \mathcal{F} \simeq \mathfrak{p}_! \Psi \Delta^* \mathcal{F}.$$
\end{lem}
\begin{proof}
We will prove $\mathfrak{p}_! \Psi \Delta_{n+1}^*  \mathcal{F} \simeq \mathfrak{p}_! \Psi \Delta_{n}^* \mathcal{F}$ for any $n \in \mathbb{Z}_{\geq 0}$. Repeat this, we deduce that $\mathfrak{p}_! \Psi \Delta_{n+d}^*  \mathcal{F} \simeq \mathfrak{p}_! \Psi \Delta_{n}^* \mathcal{F}.$ The statement of the lemma is the special case where $n=0$.

For the proof we need the stacks of shtukas with intermediate modifications $\Cht^{(1, 2)}$ and $\Cht^{(2, 1)}$.
Recall that we have a commutative diagram: 
\begin{equation*}
\xymatrixrowsep{2pc}
\xymatrixcolsep{4pc}
\xymatrix{
\Cht^{(1, 2)}  \ar[r]^{ \mathrm{Frob}_{\{1\}} }  \ar[d]^{\pi}  \ar@/_2pc/[dd]_{\mathfrak{p}^{(1, 2)}}
& \Cht^{(2, 1)}   \ar[d]^{\pi}  \ar@/^2pc/[dd]^{\mathfrak{p}^{(2, 1)}}  \\
\Cht  \ar[d]^{\mathfrak{p}}
& \Cht \ar[d]^{\mathfrak{p}}  \\
X^2  \ar[r]^{ \mathrm{Frob}_{\{1\}} } 
& X^2
}
\end{equation*}
The morphism $\pi$ is proper and small. By \ref{subsection-action-of-partial-Frob-example} we have $\pi_! \mathcal{F}^{(1, 2)} = \mathcal{F}$, $\pi_! \mathcal{F}^{(2, 1)} = \mathcal{F}$ and
\begin{equation}   \label{equation-Frob-F-2-1-equal-F-1-2}
\mathrm{Frob}_{\{1\}}^* \mathcal{F}^{(2, 1)} \simeq \mathcal{F}^{(1, 2)}.
\end{equation}

The commutative diagram
$$
\xymatrixrowsep{2pc}
\xymatrixcolsep{4pc}
\xymatrix{
X^2 \ar[r]^{\mathrm{Frob}_{\{1\}} }
& X^2 \\
\Delta_n(X) \ar[r]^{\mathrm{Frob}_{\{1\}} }  \ar@{^{(}->}[u]
& \Delta_{n+1}(X)  \ar@{^{(}->}[u]
}
$$
induces the commutative diagram
$$
\xymatrixrowsep{2pc}
\xymatrixcolsep{4pc}
\xymatrix{
\Cht^{(1, 2)}  \ar[r]^{ \mathrm{Frob}_{\{1\}} } 
& \Cht^{(2, 1)}   \\
\Cht^{(1, 2)} | _{\Delta_n} \ar[r]^{\mathrm{Frob}_{\{1\}} }  \ar@{^{(}->}[u]
& \Cht^{(2, 1)} |_{\Delta_{n+1}}   \ar@{^{(}->}[u]
}
$$
We deduce that 
\begin{equation}  \label{equation-Delta-n-Frob-equal-Frob-Delta-n+1}
\Delta_n^* \mathrm{Frob}_{\{1\}}^* \mathcal{F}^{(2, 1)} \simeq \mathrm{Frob}_{\{1\}}^* \Delta_{n+1}^* \mathcal{F}^{(2, 1)}
\end{equation}

Combining (\ref{equation-Delta-n-Frob-equal-Frob-Delta-n+1}) and (\ref{equation-Frob-F-2-1-equal-F-1-2}), we have 
$$\Delta_n^* \mathcal{F}^{(1, 2)} \simeq \mathrm{Frob}_{\{1\}}^* \Delta_{n+1}^* \mathcal{F}^{(2, 1)} $$
Since $\mathrm{Frob}_{\{1\}}$ is a homeomorpism, it commutes with nearby cycles. We have
\begin{equation} \label{equation-Psi-Delta-n-equal-Frob-Psi-Delta-n+1}
\Psi \Delta_n^* \mathcal{F}^{(1, 2)} \simeq \mathrm{Frob}_{\{1\}}^* \Psi \Delta_{n+1}^* \mathcal{F}^{(2, 1)} 
\end{equation}

We have
\begin{align*}
\mathfrak{p}_! \Psi \Delta_{n+1}^* \mathcal{F}
& \simeq \mathfrak{p}_! \Psi \Delta_{n+1}^* \pi_! \mathcal{F}^{(2, 1)} \\
& \simeq \mathfrak{p}_! \pi_! \Psi \Delta_{n+1}^* \mathcal{F}^{(2, 1)}  \\
& \simeq \mathfrak{p}^{(2, 1)}_!  \Psi \Delta_{n+1}^* \mathcal{F}^{(2, 1)}  \\
& \simeq \mathfrak{p}^{(2, 1)}_! (\mathrm{Frob}_{\{1\}})_* (\mathrm{Frob}_{\{1\}})^* \Psi \Delta_{n+1}^* \mathcal{F}^{(2, 1)}  \\
& \simeq \mathfrak{p}^{(2, 1)}_! (\mathrm{Frob}_{\{1\}})_* \Psi \Delta_n^* \mathcal{F}^{(1, 2)}  \\
& \simeq \mathfrak{p}^{(1, 2)}_! \Psi \Delta_n^* \mathcal{F}^{(1, 2)} \\
& \simeq \mathfrak{p}_! \pi_! \Psi \Delta_n^* \mathcal{F}^{(1, 2)}  \\
& \simeq \mathfrak{p}_! \Psi \Delta_n^* \mathcal{F} 
\end{align*}
The first isomorphism is because $\pi_! \mathcal{F}^{(2, 1)} = \mathcal{F}$. The second isomorphism is because $\pi$ is proper, so $\pi_! \Psi \simeq \Psi \pi_!$. The fourth isomorphism is because $(\mathrm{Frob}_{\{1\}})_* (\mathrm{Frob}_{\{1\}})^* \simeq \mathrm{Id}$. The fifth isomorphism is (\ref{equation-Psi-Delta-n-equal-Frob-Psi-Delta-n+1}). The sixth isomorphism comes from the following commutative diagram
$$
\xymatrixrowsep{2pc}
\xymatrixcolsep{4pc}
\xymatrix{
\Cht^{(1, 2)} | _{\Delta_n(x)}  \ar[r]^{ \mathrm{Frob}_{\{1\}} }  \ar[d]^{\mathfrak{p}^{(1, 2)}}
& \Cht^{(2, 1)} |_{\Delta_{n+1}(x)}  \ar[d]^{\mathfrak{p}^{(2, 1)}} \\
\Delta_n(x) \ar[r]^{\mathrm{Frob}_{\{1\}} }  
& \Delta_{n+1}(x)  
}
$$
note that $\mathrm{Frob}_{\{1\}}(x)=x$ so the lower line in this diagram is identity. Besides, since $\mathrm{Frob}_{\{1\}}$ is a homeomorphism, we have $(\mathrm{Frob}_{\{1\}})_! = (\mathrm{Frob}_{\{1\}})_*$.
The last isomorphism is because that $\pi_! \mathcal{F}^{(1, 2)} = \mathcal{F}$ and that $\pi$ is proper.
\end{proof}

Combining Lemma \ref{lem-two-legs-Psi-delta-d-eta-equal-Psi-1-Psi-2} and Lemma \ref{lem-p-Frob-d-Delta-equal-p-Delta}, we deduce $\alpha: \mathfrak{p}_! \Psi \Delta^* \mathcal{F} \xrightarrow{\sim} \mathfrak{p}_! \Psi_1 \Psi_2 \mathcal{F}$. 

\subsubsection{Nearby cycles over the base} 

\begin{lem}   \label{lem-beta-key}
Let $\mathcal{G}$ be a constant sheaf over $\overline{\eta} \times \overline{\eta}$, then 
there is a canonical morphism
\begin{equation}
\beta: \Psi \Delta^* \mathcal{G} \rightarrow  \Psi_1 \Psi_2 \mathcal{G}   
\end{equation}
which is an isomorphism.
\end{lem}
\begin{proof}
Since $\mathcal{G}$ is a constant sheaf, it is of the form $\mathcal{G} = \mathcal{A}_1 \boxtimes \mathcal{A}_2$.
Besides, over the base, $\Psi_{\overline{\eta} \rightarrow x}$ equals canonically to the fiber on $\overline{\eta}$. 
We have canonical isomorphisms
$$\Psi \Delta^* \mathcal{G} \simeq (\Delta^* \mathcal{G})_{\overline{\eta}} \simeq \mathcal{G}_{\Delta(\overline{\eta})} \simeq (\mathcal{A}_{1, \overline{\eta}}) \otimes (\mathcal{A}_{2, \overline{\eta}}) \simeq \Psi_1 \Psi_2 \mathcal{G}$$
\end{proof}

By \cite{XueSmoothness} (where we use Drinfeld's lemma), $\mathfrak{p}_! \mathcal{F}$ is constant over $\overline{\eta} \times \overline{\eta}$. By Lemma \ref{lem-beta-key}, we deduce $\beta: \Psi \Delta^* \mathfrak{p}_! \mathcal{F} \xrightarrow{\sim} \Psi_1 \Psi_2 \mathfrak{p}_! \mathcal{F}$.

\subsubsection{Commutativity}

\begin{lem}   \label{lem-alpha-key}
The following diagram is commutative 
\begin{equation*}
\xymatrix{
\mathfrak{p}_! \Psi \Delta^* \mathcal{F}   \ar[rr]^{\mathrm{can}}  \ar[d]^{\alpha}_{\simeq}
& & \Psi \Delta^* \mathfrak{p}_! \mathcal{F} \ar[d]^{\beta}_{\simeq} \\
\mathfrak{p}_! \Psi_1 \Psi_2 \mathcal{F}  \ar[r]^{\mathrm{can}}
& \Psi_1 \mathfrak{p}_!  \Psi_2 \mathcal{F} \ar[r]^{\mathrm{can}}
& \Psi_1 \Psi_2 \mathfrak{p}_! \mathcal{F} 
}
\end{equation*}
\end{lem}

\subsection{Our case: nearby cycle commutes with the partial diagonal restriction}  \label{subsection-ourcase-oneleg}

\subsubsection{Nearby cycles on shtukas}
Let $\mathfrak{p}: \Cht_{\{1, 2, 3\}} \rightarrow X^3$. Let $\Delta^{\{1, 2, 3\}}: X \rightarrow X^3$ be the diagonal morphism. Let $$\Delta^{\{1, 2\}}: X^2 \rightarrow X^3, (x, y) \mapsto (x, x, y)$$ 
$$\Delta^{\{2, 3\}}: X^2 \rightarrow X^3, (x, y) \mapsto (x, y, y)$$ be the partial diagonal morphisms. For any $d \in \mathbb{N}$, let 
$\Delta^{\{1, 2, 3\}}_d := \mathrm{Frob}_{\{1\}}^{2d} \mathrm{Frob}_{\{2\}}^{d} \Delta^{\{1, 2, 3\}}$,  
$\Delta^{\{1, 2\}}_d := \mathrm{Frob}_{\{1\}}^{d} \Delta^{\{1, 2\}}$, 
$\Delta^{\{2, 3\}}_d := \mathrm{Frob}_{\{2\}}^{d} \Delta^{\{2, 3\}}$. 
Let $\mathcal{F} = \mathcal{S}_{\{1, 2, 3\}, W_1 \boxtimes W_2 \boxtimes W_3} \otimes \mathcal{L}^{nx}$.

\begin{lem}   \label{lem-psi-Delta-psi-1-2-3}
For $d \gg 0$, we have the following canonical morphisms which are isomorphisms and such that the following diagram is commutative
\begin{equation}   \label{equation-Psi-Delta-1-2-3-triangle-diagram}
\xymatrix{
& \Psi_1 \Psi_2 (\Delta^{\{2, 3\}}_d)^* \mathcal{F} \ar[d]^{\simeq} \\
\Psi (\Delta^{\{1, 2, 3\}}_d)^* \mathcal{F}   \ar[r]^{\simeq}  \ar[ru]^{\simeq}  \ar[rd]^{\simeq}
& \Psi_1 \Psi_2 \Psi_3 \mathcal{F}  \\
& \Psi_2 (\Delta^{\{1, 2\}}_d)^* \Psi_3 \mathcal{F}  \ar[u]^{\simeq}
}
\end{equation}
\end{lem}
\begin{proof}
The proof is similar to the proof of Lemma \ref{lem-two-legs-Psi-delta-d-eta-equal-Psi-1-Psi-2}.
\end{proof}

\begin{constr}  \label{constr-morphism-alpha}
Similar to Lemma \ref{lem-p-Frob-d-Delta-equal-p-Delta}, when we apply $\mathfrak{p}_!$ to (\ref{equation-Psi-Delta-1-2-3-triangle-diagram}), we can remove the index $d$.
In particular, we have composition of morphisms
\begin{equation}   \label{equation-construction-alpha}
\alpha: \mathfrak{p}_!  \Psi_2 (\Delta^{\{2, 3\}})^* \mathcal{F}_{\{1, 2, 3\}} \xrightarrow{\sim} \mathfrak{p}_! \Psi_2 (\Delta_d^{\{2, 3\}})^* \mathcal{F}_{\{1, 2, 3\}}  \xrightarrow{\sim} \mathfrak{p}_!  \Psi_2 \Psi_3 \mathcal{F}_{\{1, 2, 3\}} 
\end{equation}
\end{constr}

\subsubsection{Nearby cycles over the base}  \label{subsectionDrinfeldlemma-oneleg}

\begin{prop}   \label{prop-p-F-constant}
(1) $\mathfrak{p}_!  \mathcal{F}_{\{1, 2, 3\}}$ is constant over $\overline{\eta} \times \overline{\eta} \times \overline{\eta}$.

(2) $\mathfrak{p}_!   \Psi_3 \mathcal{F}_{\{1, 2, 3\}}$ is constant over $\overline{\eta} \times \overline{\eta} \times x$.
\end{prop}
\begin{proof}
In fact, $\mathfrak{p}_! \mathcal{F}_{\{1, 2, 3\}}$ and $\mathfrak{p}_!   \Psi_3 \mathcal{F}_{\{1, 2, 3\}}$ are equipped with an action of partial Frobenius morphisms and have the Eichler-Shimura relations. We apply Drinfeld's lemma as in \cite{XueSmoothness}. We send to \ref{subsectionFusiononCurve} for details.
\end{proof}

\begin{constr}   \label{constr-morphism-beta}
Applying Lemma \ref{lem-beta-key} to $\mathcal{G} = \mathfrak{p}_! \Psi_3 \mathcal{F}_{\{1, 2, 3\}}$, whose condition is satisfied by Proposition \ref{prop-p-F-constant},
we construct a morphism
\begin{equation}   \label{equation-construction-beta}
\beta: \Psi (\Delta^{\{1, 2\}})^* \mathfrak{p}_!  \Psi_3 \mathcal{F}_{\{1, 2, 3\}} \xrightarrow{\sim} 
\Psi_1 \Psi_2 \mathfrak{p}_!   \Psi_3 \mathcal{F}_{\{1, 2, 3\}} 
\end{equation}
\end{constr}

\subsubsection{Commutativity}

\begin{lem}    \label{lem-p-Psi-1-2-3-big-commutative-diagram}
We have the following canonical morphisms such that the following diagram is commutative
$${  \resizebox{15cm}{!}{ 
\!\!\!\!\!\!\!
\xymatrix{
& \mathfrak{p}_! \Psi_1 \Psi_2 (\Delta^{23})^* \mathcal{F} \ar[d]^{\simeq}_{\alpha}   \ar[r]^{\mathrm{can}}
& \Psi_1 \mathfrak{p}_! \Psi_2 (\Delta^{23})^* \mathcal{F} \ar[d]^{\simeq}_{\alpha}  \ar[rr]^{\mathrm{can}}  \ar@{}[rrd]|-{(b)} 
& & \Psi_1 \Psi_2 (\Delta^{23})^* \mathfrak{p}_! \mathcal{F} \ar[d]^{\simeq}_{\beta} \ar[rd]^{\simeq}_{\beta}  \\
\mathfrak{p}_! \Psi (\Delta^{123})^* \mathcal{F}   \ar[r]^{\simeq}  \ar[ru]^{\simeq}_{\alpha}  \ar[rd]_{\simeq}^{\alpha}
& \mathfrak{p}_! \Psi_1 \Psi_2 \Psi_3 \mathcal{F}  \ar[r]^{\mathrm{can}}
&  \Psi_1 \mathfrak{p}_! \Psi_2 \Psi_3 \mathcal{F}  \ar[r]^{\mathrm{can}}
&  \Psi_1 \Psi_2 \mathfrak{p}_! \Psi_3 \mathcal{F}  \ar[r]^{\mathrm{can}}
&  \Psi_1 \Psi_2 \Psi_3 \mathfrak{p}_! \mathcal{F}   \ar[r]^{\simeq}
& \Psi (\Delta^{123})^* \mathfrak{p}_! \mathcal{F}   \\
& \mathfrak{p}_! \Psi_2 (\Delta^{12})^* \Psi_3 \mathcal{F}  \ar[u]_{\simeq}^{\alpha}  \ar[rr]^{\mathrm{can}}   \ar@{}[rru]|-{(a)} 
& & \Psi_2 (\Delta^{12})^* \mathfrak{p}_! \Psi_3  \mathcal{F}  \ar[u]_{\simeq}^{\beta}  \ar[r]^{\mathrm{can}}
& \Psi_2 (\Delta^{12})^* \Psi_3 \mathfrak{p}_!  \mathcal{F}  \ar[u]_{\simeq}^{\beta}   \ar[ru]^{\simeq}_{\beta}
}
}
}
$$
where $\alpha$ are morphisms coming from nearby cycles over stacks of shtukas, $\beta$ are morphisms coming from nearby cycles over the base.
\end{lem}
\begin{proof}
Similar to Lemma \ref{lem-alpha-key}. We combine Lemma \ref{lem-psi-Delta-psi-1-2-3} and Proposition \ref{prop-p-F-constant}.
\end{proof}

\subsection{``Zorro lemma" argument}  \label{subsectionZorro-oneleg}

\begin{lem} (Zorro lemma)  \label{lem-Zorro-one-leg}
The composition $$\mathcal{F}_{\{1\}, W} \otimes \Lambda \xrightarrow{ \mathcal{C}^{\sharp}_{\{2, 3\} }  } (\Delta^{\{1, 2, 3\}})^*\mathcal{F}_{\{1, 2, 3\}, W \boxtimes W^* \boxtimes W} \xrightarrow{ \mathcal{C}^{\flat}_{\{1, 2\}} } \Lambda \otimes \mathcal{F}_{\{3\}, W}$$ 
is the identity.
\end{lem}
\begin{proof}
In fact, by the fusion property of Satake sheaves, $(\Delta^{\{1, 2, 3\}})^*\mathcal{F}_{\{1, 2, 3\}, W \boxtimes W^* \boxtimes W} = \mathcal{F}_{\{1\}, W \otimes W^* \otimes W}$. We know that the composition of morphisms of vector spaces $W \xrightarrow{\mathrm{Id} \otimes \delta} W \otimes W^* \otimes W \xrightarrow{\mathrm{ev} \otimes \mathrm{Id}} W$ is the identity. The lemma then follows from the functoriality of $\mathcal{F}_{I, W}$ on $W$.
\end{proof}

\quad

\begin{lem} 
The composition $\gamma \circ \mathrm{can}$ is an isomorphism ($\gamma$ is constructed in  \ref{subsectionInversemap}).
\end{lem}

\begin{proof}
Step (1): The following diagram is commutative, where the composition of  the right vertical morphisms is $\gamma$:
\begin{equation}   \label{diagram-injectivity}
\xymatrix{
\mathfrak{p}_! \Psi_1 \mathcal{F}_{\{1\}, W}  \otimes \Lambda \ar[r]^{\mathrm{can}}   \ar[d]^{\mathcal{C}^{\sharp}_{\{2, 3\}}} 
& \Psi_1 \mathfrak{p}_! \mathcal{F}_{\{1\}, W}  \otimes \Lambda  \ar[d]^{\mathcal{C}^{\sharp}_{\{2, 3\}}} \\
\mathfrak{p}_! \Psi_1 \Psi_2 (\Delta^{23})^* \mathcal{F}_{\{1, 2, 3\}, W \boxtimes W^* \boxtimes W}   \ar[r]^{\mathrm{can}}  \ar[d]^{\alpha}
& \Psi_1 \mathfrak{p}_!  \Psi_2 (\Delta^{23})^* \mathcal{F}_{\{1, 2, 3\}, W \boxtimes W^* \boxtimes W}  \ar[d]^{\alpha} \\
\mathfrak{p}_! \Psi_1 \Psi_2 \Psi_3 \mathcal{F}_{\{1, 2, 3\}, W \boxtimes W^* \boxtimes W}   \ar[r]^{\mathrm{can}}  \ar[rd]_{\mathrm{can}}  \ar[dd]^{\alpha^{-1}}
& \Psi_1 \mathfrak{p}_!  \Psi_2 \Psi_3 \mathcal{F}_{\{1, 2, 3\}, W \boxtimes W^* \boxtimes W} \ar[d]^{\mathrm{can}} \\
& \Psi_1 \Psi_2 \mathfrak{p}_!   \Psi_3 \mathcal{F}_{\{1, 2, 3\}, W \boxtimes W^* \boxtimes W}  \ar[d]^{\beta^{-1}} \\
\mathfrak{p}_! \Psi_2 (\Delta^{12})^* \Psi_3 \mathcal{F}_{\{1, 2, 3\}, W \boxtimes W^* \boxtimes W}  \ar[r]^{\mathrm{can}}  \ar[d]^{\mathcal{C}^{\flat}_{\{1, 2\}}}   \ar@{}[ru]^{(a)}   
& \Psi_2 (\Delta^{12})^* \mathfrak{p}_!  \Psi_3 \mathcal{F}_{\{1, 2, 3\}, W \boxtimes W^* \boxtimes W}  \ar[d]^{\mathcal{C}^{\flat}_{\{1, 2\}}}  \\
\Lambda \otimes \mathfrak{p}_! \Psi_3 \mathcal{F}_{\{3\}, W}  \ar[r]^{\text{Id}}
& \Lambda \otimes \mathfrak{p}_! \Psi_3 \mathcal{F}_{\{3\}, W}
}
\end{equation}

Note that the first square is:
$$
\xymatrix{
\mathfrak{p}_! \Psi_1 \mathcal{F}_{\{1\}, W}  \otimes \Lambda \ar[r]^{\mathrm{can}}   \ar[d]^{ \simeq }    \ar@/_3pc/[ddd]_{\mathcal{C}^{\sharp}_{\{2, 3\}}} 
& \Psi_1 \mathfrak{p}_! \mathcal{F}_{\{1\}, W}  \otimes \Lambda  \ar[d]^{ \simeq }  \ar@/^3pc/[ddd]^{\mathcal{C}^{\sharp}_{\{2, 3\}}} \\
\mathfrak{p}_! \Psi_1 \Psi_2 \mathcal{F}_{\{1, 2\}, W \boxtimes \mathbf{1}}  \ar[r]^{\mathrm{can}}   \ar[d]^{\delta}
& \Psi_1 \mathfrak{p}_! \Psi_2 \mathcal{F}_{\{1, 2\}, W \boxtimes \mathbf{1}}   \ar[d]^{\delta}  \\
\mathfrak{p}_! \Psi_1 \Psi_2 \mathcal{F}_{\{1, 2\}, W \boxtimes (W^* \otimes W)}  \ar[r]^{\mathrm{can}}   \ar[d]^{\simeq}
& \Psi_1 \mathfrak{p}_! \Psi_2 \mathcal{F}_{\{1, 2\}, W \boxtimes (W^* \otimes W)}   \ar[d]^{\simeq}  \\
\mathfrak{p}_! \Psi_1 \Psi_2 (\Delta^{23})^* \mathcal{F}_{\{1, 2, 3\}, W \boxtimes W^* \boxtimes W}   \ar[r]^{\mathrm{can}} 
& \Psi_1 \mathfrak{p}_!  \Psi_2 (\Delta^{23})^* \mathcal{F}_{\{1, 2, 3\}, W \boxtimes W^* \boxtimes W} 
}
$$
The last square is:
$$
\xymatrix{
\mathfrak{p}_! \Psi_2 (\Delta^{12})^* \Psi_3 \mathcal{F}_{\{1, 2, 3\}, W \boxtimes W^* \boxtimes W}  \ar[r]^{\mathrm{can}}  \ar[d]^{\simeq}   \ar@/_3pc/[ddd]_{\mathcal{C}^{\flat}_{\{1, 2\}}} 
& \Psi_2 (\Delta^{12})^* \mathfrak{p}_!  \Psi_3 \mathcal{F}_{\{1, 2, 3\}, W \boxtimes W^* \boxtimes W}  \ar[d]^{\simeq} \ar@/^3pc/[ddd]^{\mathcal{C}^{\flat}_{\{1, 2\}}}  \\
\mathfrak{p}_! \Psi_2 \Psi_3 \mathcal{F}_{\{2, 3\}, (W \otimes W^*) \boxtimes W}  \ar[r]^{\mathrm{can}}  \ar[d]^{\mathrm{ev}}
& \Psi_2 \mathfrak{p}_!  \Psi_3 \mathcal{F}_{\{2, 3\}, (W \otimes W^*) \boxtimes W}  \ar[d]^{\mathrm{ev}}  \\
\mathfrak{p}_! \Psi_2 \Psi_3 \mathcal{F}_{\{2, 3\}, \mathbf{1} \boxtimes W}  \ar[r]^{\mathrm{can}}  \ar[d]^{\simeq}
& \Psi_2 \mathfrak{p}_!  \Psi_3 \mathcal{F}_{\{2, 3\}, \mathbf{1} \boxtimes W}  \ar[d]^{\simeq}  \\
\Lambda \otimes \mathfrak{p}_! \Psi_3 \mathcal{F}_{\{3\}, W}  \ar[r]^{\text{Id}}
& \Lambda \otimes \mathfrak{p}_! \Psi_3 \mathcal{F}_{\{3\}, W}
}
$$

The only commutativity needs to prove is (a), this follows from Lemma \ref{lem-p-Psi-1-2-3-big-commutative-diagram}. The other squares are commutative because the canonical morphism "can" is functorial, so it commutes with morphisms of sheaves.

Step (2): The composition of the left vertical line of (\ref{diagram-injectivity}) is the identity.

In fact, the following diagram is commutative. This follows from Lemma \ref{lem-p-Psi-1-2-3-big-commutative-diagram} and that the morphism $\alpha$ is functorial, so it commutes with morphisms of sheaves.
\begin{equation}   \label{diagram-injectivity-Zorro}
\xymatrix{
\mathfrak{p}_! \Psi_1 \mathcal{F}_{\{1\}, W}  \otimes \Lambda \ar[r]^{\mathrm{Id}}     \ar[dd]^{\mathcal{C}^{\sharp}_{\{2, 3\}}} 
& \mathfrak{p}_! \Psi_1 \mathcal{F}_{\{1\}, W}  \otimes \Lambda  \ar[d]^{\mathcal{C}^{\sharp}_{\{2, 3\}}}  \\
& \mathfrak{p}_! \Psi_1 \Psi_2 (\Delta^{23})^* \mathcal{F}_{\{1, 2, 3\}, W \boxtimes W^* \boxtimes W}    \ar[d]^{\simeq}_{\alpha} \\
\mathfrak{p}_! \Psi (\Delta^{1,2,3})^* \mathcal{F}_{\{1, 2, 3\}, W \boxtimes W^* \boxtimes W}   \ar[r]^{\simeq}_{\alpha}  \ar[ru]^{\simeq}_{\alpha}   \ar[rd]^{\simeq}_{\alpha}   \ar[dd]^{\mathcal{C}^{\flat}_{\{1, 2\}}} 
& \mathfrak{p}_! \Psi_1 \Psi_2 \Psi_3 \mathcal{F}_{\{1, 2, 3\}, W \boxtimes W^* \boxtimes W} \ar[d]^{\simeq}_{\alpha^{-1}}  \\
& \mathfrak{p}_! \Psi_2 (\Delta^{12})^* \Psi_3 \mathcal{F}_{\{1, 2, 3\}, W \boxtimes W^* \boxtimes W}  \ar[d]^{\mathcal{C}^{\flat}_{\{1, 2\}}}   \\
\Lambda \otimes \mathfrak{p}_! \Psi_3 \mathcal{F}_{\{3\}}  \ar[r]^{\text{Id}}
& \Lambda \otimes \mathfrak{p}_! \Psi_3 \mathcal{F}_{\{3\}}
}
\end{equation}
The left vertical line of (\ref{diagram-injectivity}) is the right vertical line of (\ref{diagram-injectivity-Zorro}). By Lemma \ref{lem-Zorro-one-leg}, the composition of the left vertical line of (\ref{diagram-injectivity-Zorro}) is the identity.
\end{proof}

\begin{lem}
The composition $\mathrm{can} \circ \gamma$ is an isomorphism.
\end{lem}

\begin{proof}
Step (1): The following diagram is commutative, where the composition of the left vertical line is $\gamma$:
\begin{equation}   \label{diagram-surjectivity}
\xymatrix{
\Psi_1 \mathfrak{p}_! \mathcal{F}_{\{1\}, W}  \otimes \Lambda  \ar[d]^{\mathcal{C}^{\sharp}_{\{2, 3\}}}   \ar[r]^{\mathrm{Id}}
& \Psi_1 \mathfrak{p}_! \mathcal{F}_{\{1\}, W}  \otimes \Lambda \ar[d]^{\mathcal{C}^{\sharp}_{\{2, 3\}}}  \\
\Psi_1 \mathfrak{p}_!  \Psi_2 (\Delta^{\{2, 3\}})^* \mathcal{F}_{\{1, 2, 3\}, W \boxtimes W^* \boxtimes W}  \ar[r]^{\mathrm{can}}   \ar[d]^{\alpha}   \ar@{}[rd]^{(b)} 
&  \Psi_1 \Psi_2 (\Delta^{\{2, 3\}})^* \mathfrak{p}_!   \mathcal{F}_{\{1, 2, 3\}, W \boxtimes W^* \boxtimes W}  \ar[d]^{\beta} \\
\Psi_1 \mathfrak{p}_!  \Psi_2 \Psi_3 \mathcal{F}_{\{1, 2, 3\}, W \boxtimes W^* \boxtimes W} \ar[d]^{\mathrm{can}}   \ar[r]^{\mathrm{can}}   
& \Psi_1  \Psi_2 \Psi_3 \mathfrak{p}_! \mathcal{F}_{\{1, 2, 3\}, W \boxtimes W^* \boxtimes W}  \ar[dd]^{\beta^{-1}} \\
\Psi_1  \Psi_2 \mathfrak{p}_! \Psi_3 \mathcal{F}_{\{1, 2, 3\}, W \boxtimes W^* \boxtimes W} \ar[ru]^{\mathrm{can}} \ar[d]^{\beta^{-1}}  \\
\Psi_2 (\Delta^{\{1, 2\}})^* \mathfrak{p}_!  \Psi_3 \mathcal{F}_{\{1, 2, 3\}, W \boxtimes W^* \boxtimes W}  \ar[d]^{\mathcal{C}^{\flat}_{\{1, 2\}}}   \ar[r]^{\mathrm{can}}
& \Psi_2 (\Delta^{\{1, 2\}})^* \Psi_3 \mathfrak{p}_!  \mathcal{F}_{\{1, 2, 3\}, W \boxtimes W^* \boxtimes W}  \ar[d]^{\mathcal{C}^{\flat}_{\{1, 2\}}} \\
\Lambda \otimes \mathfrak{p}_! \Psi_3 \mathcal{F}_{\{3\}, W}   \ar[r]^{\mathrm{can}}
& \Lambda \otimes \Psi_3 \mathfrak{p}_! \mathcal{F}_{\{3\}, W} 
}
\end{equation}
The only commutativity needs to prove is (b), this follows from Lemma \ref{lem-p-Psi-1-2-3-big-commutative-diagram}. The other squares are commutative because the canonical morphism "can" is functorial, so it commutes with morphisms of sheaves.

Step (2): The composition of the right vertical line of (\ref{diagram-surjectivity}) is the identity.

In fact, the following diagram is commutative. This follows from Lemma \ref{lem-p-Psi-1-2-3-big-commutative-diagram} and that the morphism $\beta$ is functorial, so it commutes with morphisms of sheaves.
\begin{equation}   \label{diagram-surjectivity-Zorro}
\xymatrix{
\Psi_1 \mathfrak{p}_! \mathcal{F}_{\{1\}, W}  \otimes \Lambda  \ar[d]^{\mathcal{C}^{\sharp}_{\{2, 3\}}}   \ar[r]^{\mathrm{Id}}
& \Psi_1 \mathfrak{p}_! \mathcal{F}_{\{1\}, W}  \otimes \Lambda  \ar[dd]^{\mathcal{C}^{\sharp}_{\{2, 3\}}}   \\
\Psi_1 \Psi_2 (\Delta^{\{2, 3\}})^* \mathfrak{p}_!   \mathcal{F}_{\{1, 2, 3\}, W \boxtimes W^* \boxtimes W}  \ar[d]^{\beta}  \ar[rd]^{\simeq}_{\beta}
& \\
\Psi_1  \Psi_2 \Psi_3 \mathfrak{p}_! \mathcal{F}_{\{1, 2, 3\}, W \boxtimes W^* \boxtimes W}  \ar[d]^{\beta^{-1}}   \ar[r]^{\simeq}_{\beta}
& \Psi (\Delta^{123})^* \mathfrak{p}_! \mathcal{F}_{\{1, 2, 3\}, W \boxtimes W^* \boxtimes W} \ar[dd]^{\mathcal{C}^{\flat}_{\{1, 2\}}}  \\
\Psi_2 (\Delta^{\{1, 2\}})^* \Psi_3 \mathfrak{p}_!  \mathcal{F}_{\{1, 2, 3\}, W \boxtimes W^* \boxtimes W}   \ar[d]^{\mathcal{C}^{\flat}_{\{1, 2\}}}  \ar[ru]^{\simeq}_{\beta}
&  \\
\Lambda \otimes \Psi_3 \mathfrak{p}_! \mathcal{F}_{\{3\}, W}  \ar[r]^{\mathrm{Id}} 
& \Lambda \otimes \Psi_3 \mathfrak{p}_! \mathcal{F}_{\{3\}, W} 
}
\end{equation}
The right vertical line of (\ref{diagram-surjectivity}) is the left vertical line of (\ref{diagram-surjectivity-Zorro}). By Lemma \ref{lem-Zorro-one-leg}, the composition of the right vertical line of (\ref{diagram-surjectivity-Zorro}) is the identity.
\end{proof}

\section{Fusion properties for nearby cycles}\label{sectionFusion}

In this section, we prepare some technical results needed for the construction of the inverse map. In \ref{subsectionFusiononCurve} we generalize the construction of morphism $\beta$. In \ref{subsectionFusiononshtukas} we generalize the construction of morphism $\alpha$.

\subsection{Reminder on some properties of stacks of shtukas needed later}  \label{subsection-ChtR-intermediate}

The statements in main theorems involve only stacks of shtukas without intermediate modifications. However in the proof we need to consider stack of shtukas with intermediate modifications, for two reasons:
the twisted product structure and the action of partial Frobenius morphisms.

\subsubsection{Relations}
Let $I=I_1 \cup I_2$ and $W = W_1 \boxtimes W_2$. We denote by $\Cht_{I, W}^{(I_1, I_2)}$ (resp. $\ChtR_{I, W}^{(I_1, I_2)}$) the stack of global shtukas (resp. the stack of restricted shtukas) with intermediate modifications. We send to \cite{Lafforgue} and \cite{GenestierLafforgue} for the definition. We have the Cartesian diagram over $X^I$:
\begin{equation}   \label{equation-Cht-ChtR-intermediate}
\xymatrixrowsep{1.5pc}
\xymatrixcolsep{5pc}
\xymatrix{
\Cht_{I, W}^{(I_1, I_2)}   \ar[r]^{\mathcal{R}^{nx}}  \ar[d]^{\pi}
& \ChtR_{I, W}^{nx, (I_1, I_2)} \ar[d]^{\pi}  \ar[r]
& [L^+_I G \backslash \Gr_{I, W}^{(I_1, I_2)}] \ar[d]^{\pi} \\
\Cht_{I, W} \ar[r]^{\mathcal{R}^{nx}}
& \ChtR_{I, W}^{nx}  \ar[r]
& [L^+_I G \backslash \Gr_{I, W} ]
}
\end{equation}
The horizontal morphisms are smooth. 
The vertical morphisms are proper, small, and are isomorphism over the open subscheme of $X^I$ where the legs are two are two distinct. By the geometric Satake equivalence (\protect{\cite[Théorème 1.17 b)]{Lafforgue}}), we have canonical isomorphism $$\pi_! \mathcal{S}_{I, W_1 \boxtimes W_2}^{(I_1, I_2)} = \mathcal{S}_{I, W_1 \boxtimes W_2}.$$ It induces canonical isomorphism $$\pi_! \mathcal{F}_{I, W_1 \boxtimes W_2}^{(I_1, I_2)} = \mathcal{F}_{I, W_1 \boxtimes W_2}.$$
%where $\mathcal{F}_{I, W_1 \boxtimes W_2}$ is defined in \ref{subsection-def-F-I-W} and $\mathcal{F}_{I, W_1 \boxtimes W_2}^{(I_1, I_2)}$ is defined similarly.

\subsubsection{Product structure on restricted shtukas}
\label{subsubsectionProductstructure}

We first recall the notion of twisted product. 
\begin{defi}
\begin{enumerate}
\item Let $K$ be an algebraic group and $\mathcal{X}, \mathcal{Y}$ be two schemes or stacks, we assume that we are given the data of $\mcE \to \mathcal{X}$ a $K$-torsor on $\mathcal{X}$ and an action of $K$ on $\mathcal{Y}$. The twisted product $\mathcal{X} \tilde{\times} \mathcal{Y}$ is the stack $\mcE \times \mathcal{Y}/K$ where $K$ acts diagonally on $\mcE \times \mathcal{Y}$. Usually when we write a twisted product $\mathcal{X} \tilde{\times} \mathcal{Y}$, the $K$-torsor is implicit. 
\item Let $\mcA$ be a sheaf on $\mathcal{X}$ and $\mcB$ be a $K$-equivariant sheaf on $\mcB$, the twisted exterior product $\mcA \tilde{\boxtimes} \mcB$ is the unique sheaf on $\mathcal{X} \tilde{\times} \mathcal{Y}$ whose pullback to $\mcE \times \mathcal{Y}$ is $\mcA \boxtimes \mcB$. 
\end{enumerate}
\end{defi}

\begin{lem}\label{lemTwistedProductStructure}
Let $I_1, I_2$ be two finite sets, $W_i \in \Rep_{\Lambda} {^L}G^{I_i}$ two representations and $(n_i)_{i \in I_1 \cup I_2}$ a collection of integers such that $\ChtR^{nx, (I_1, I_2)}_{I_1 \cup I_2, W_1 \boxtimes W_2, (n_i)}$ is well defined. This is a stack over $X^{I_1 \cup I_2}$. 
\begin{enumerate}
\item Over $X^{I_1} \times \check{X}^{I_2}$ this stack splits as twisted product, that is $$\ChtR^{nx, (I_1, I_2)}_{I_1 \cup I_2, W_1 \boxtimes W_2, (n_i)}= \ChtR^{nx}_{I_1, W_1, (n_i)_{i \in I_1}} \; \tilde{\times} \; G_{\sum_{i \in I_2} n_iy_i} \backslash \Gr_{I_2, W_2}.$$ 
\item Given a sheaf $\mcA \in \ChtR^{nx}_{I_1, W_1, (n_i)_{i \in I_1}}$, the pullback to $\ChtR^{nx}_{I_1, W_1, (n_i)_{i \in I_1}} \times \Gr_{I_2, W_2}$ of the sheaf $\mcA \tilde{\boxtimes} \mcS_{I_2, W_2}$ is $\mcA {\boxtimes} \mcS_{I_2, W_2}$.
\end{enumerate}
\end{lem}

\begin{proof}
$(1).$ We will use the notations of \cite[Definition 2.10]{GenestierLafforgue} Let $((y_i)_{i \in I_1\cup I_2}, \mcG, z, \theta)$ be an $S$-point of $\ChtR^{nx, (I_1, I_2)}_{I_1 \cup I_2, W_1 \boxtimes W_2, (n_i)}$ with $(y_i) \in (X^{I_1} \times \check{X}^{I_2})(S)$. And recall that $z$ is a point of $\Gr_{I_1 \cup I_2, W_1 \boxtimes W_2}^{(I_1, I_2), (y_i)} \times_{G_{\sum n_iy_i}} \mcG_{\Gamma_{\sum n_iy_i}}$, that is it parametrizes a sequence of modification of $G$-torsors on $X$
$$\mcE_1 \to \mcE_2 \to \mcE^0$$
happening at the legs indexed by $I_1$ and $I_2$ respectively, where $\mcE^0$ denotes the trivial $G$-torsor and an identification of $\mcE_1$ on $\sum n_ix_i$ with $\mcG$. 
Let $\mcE_{\univ} \to \ChtR^{nx, (I_1, I_2)}_{I_1 \cup I_2, W_1 \boxtimes W_2, (n_i)}$ be the $G_{\sum_{i \in I_2} n_iy_i}$-torsor parametrizing trivializations of $\mcE_2$ on $\sum_{i \in I_2} n_iy_i$. 
Given $((y_i)_{i \in I_1\cup I_2}, \mcG, z, \theta, \psi)$ an $S$-point of $\mcE_{\univ}$ where $\psi$ denotes the trivialization of $\mcE_2$, denote by $z_1$ the point of $\Gr_{I_1 , W_1}^{(I_1), (y_i)_{i \in I_1}} \times_{G_{\sum_{i \in I_1} n_iy_i}} \mcG_{\Gamma_{\sum_{I_1} n_iy_i}}$ obtained by composing the modification $\mcE_1 \to \mcE_2$ with the trivialization $\psi$. Using the notations $a_z(\mcG)$ of \cite[Notation 2.8]{GenestierLafforgue}, we note that there is a canonical isomorphism of $G$-torsor on $nx$
$$a_z(\mcG) = a_{z_1}(\mcG)$$
since the second modifications happens away from $x$. There is therefore a well defined map
$$\mcE_{\univ} \to \ChtR^{nx}_{I_1, W_1, (n_i)}$$
obtained by sending $((y_i)_{i \in I_1\cup I_2}, \mcG, z, \theta, \psi)$ to $((y_i)_{i \in I_1}, \mcG_{|nx + \sum_{i \in I_1} n_iy_i}, z_1, \theta_1)$ where $\theta_1$ is the map 
$${^\sigma}a_{z_1}(\mcG) = {^\sigma}a_z(\mcG) \xrightarrow{\theta} \mcG_{|nx}.$$
On the other hand the trivialization $\psi$ defines a map $\mcE_{\univ} \to \Gr_{I_2, W_2}$ over $\check{X}^{I_2}$ compatible with the $G_{\sum_{i \in I_2} n_iy_i}$-action. Finally it is clear that that map $\mcE_{\univ} \to \Gr_{I_2, W_2} \times \ChtR^{nx}_{I_1, W_1, (n_i)}$ is an isomorphism. 
$(2).$ The statement about the Satake sheaves is clear. 
\end{proof}

\begin{rque}   \label{rem-ChtR-product}
In the context of the above lemma, over $x^{I_1} \times \check{X}^{I_2}$, the twisted product is a product. Besides, over $x^{I_1} \times \check{X}^{I_2}$ the morphisms $\pi$ are isomorphisms. We have
%\label{equation-product-over-X-times-x}
\begin{equation*}  
\ChtR^{nx}_{I, W, (n_i)} |_{x^{I_1} \times \check{X}^{I_2}} = \ChtR^{nx, (I_1, I_2)}_{I, W, (n_i)} |_{x^{I_1} \times \check{X}^{I_2}}= \ChtR^{nx}_{I_1, W_1, (n_i)_{i \in I_1}} |_{x^{I_1}} \times [G_{\sum_{i \in I_2} n_iy_i} \backslash \Gr_{I_2, W_2}] |_{\check{X}^{I_2}}.
\end{equation*}
Morphism (\ref{equation-R-nx-intro}) is the composition of $\Cht_{I \cup \{0\}, W \boxtimes W_0} \rightarrow \ChtR_{I \cup \{0\}, W \boxtimes W_0, (n_i)}^{nx}$ and the above morphism, for $I_1=\{0\}, W_1 = W_0$ and $I_2=I, W_2 = W$.
\end{rque}

\subsubsection{Removing the level structures outside of $x$}

\begin{lem}\label{lemLevelStructure}
For all modifications $\widetilde{S} \to X^I, s \in \widetilde{S}$ above $x$ and specialization maps $\overline{\eta}_I \to s$ there is a canonical isomorphism 
$$\Psi_{\overline{\eta}_I \to s} \mcF_{I,W, N^x} = (\mcL_{N^x})_{(\Cht)_{|s}} \otimes \Psi_{\overline{\eta}_I \to s} \mcF_{I,W, \emptyset}.$$
\end{lem}

\begin{proof}
Since $\mcF_{I,W, N^x} = \mcL_{N^x} \otimes \mcF_{I,W, \emptyset}$ and $\mcL_{N^x}$ is locally constant in a neighborhood of $x$ hence of $s$, the canonical map 
$$\Psi_{\overline{\eta}_I \to s}\mcL_{N^x} \otimes \Psi_{\overline{\eta}_I \to s}\mcF_{I,W, \emptyset} \to \Psi_{\overline{\eta}_I \to s} \mcF_{I,W, N^x}$$
coming from the lax monoidality of $\Psi_{\overline{\eta}_I \to s}$ is an isomorphism. Furthermore $\Psi_{\overline{\eta}_I \to s}\mcL_{N^x}$ is isomorphic to $(\mcL_{N^x})_{(\Cht)_{|s}}$. 
\end{proof}

By lemma \ref{lemLevelStructure}, we can assume that $N^x = \emptyset$ and therefore the sheaf $\mcF_{I, W}$ is then pullback from the corresponding on the stack $\ChtR_{I \cup J \cup \{0\}, (n_i)}^{nx}$ along the restriction map.

\subsubsection{Action of partial Frobenius morphisms}
 \label{subsection-action-of-partial-Frob-example}

\begin{defi}
Let $I$ be a finite set, for $d = (d_i) \in \N^I$ we denote by $\Frob^d : X^I \to X^I$ the morphism $(x_i) \mapsto (\Frob^{d_i}(x_i))$ where $\Frob : X \to X$ is the $q$-power $\Fqb$-linear Frobenius. We call these morphism partial Frobenius, classically the partial Frobenius morphisms are the morphisms $\Frob^d$ when $d = (0, \dots, 0,1,0, \dots, 0)$ has only one nonzero entry which is a $1$ at the $i$-th item, and is classically denoted by $\Frob_i$. 
\end{defi}

Recall that the action of the partial Frobenius morphisms on a complex $\mathcal{G}$ over $X^I$ is the following data: for every $i \in I$, a morphism $\Frob_i: \Frob_i^* \mathcal{G} \rightarrow \mathcal{G}$ commuting between them such that the composition is the total Frobenius morphism on $\mathcal{G}$.

For $\mathcal{F}_{I, W} = \mathcal{R}^* (\mathcal{S}_{I, W} \boxtimes \mathcal{A})$ over $\Cht_{I \cup \{0\}, W \boxtimes W_0} |_{\check{X}^I \times x}$, the proper direct image $\mathfrak{p}_! \mathcal{F}_{I, W}$ is equipped with an action of partial Frobenius morphisms.

The construction of this action of partial Frobenius morphisms is the same as in \protect{\cite[Section 3]{Lafforgue}} or \protect{\cite[Section 7.1]{XueIntegral}}, except that we replace 
\begin{equation}  \label{equation-epsilon}
\widetilde{\epsilon}: \Cht_{I, W}^{(I_1, \cdots, I_k)} \rightarrow \prod_{j=1}^k [L^+_{I_j}G \backslash \Gr_{I_j, W_j}]
\end{equation}
by 
\begin{equation} \label{equation-R}
\mathcal{R}: \Cht_{I \cup \{0\}, W \boxtimes W_0}^{(I_1, \cdots, I_k, \{0\})} |_{\check{X}^I \times x} \rightarrow \prod_{j=1}^k [L^+_{I_j}G \backslash \Gr_{I_j, W_j}] |_{\check{X}^I} \times \ChtR^{nx}_{\{0\}, W_0} |_{x}
\end{equation}
and replace $\mathcal{F}_{I, W}^{(I_1, \cdots, I_k)} = \widetilde{\epsilon}^*(\boxtimes_j \mathcal{S}_{I_j, W_j})$ by $\mathcal{F}_{I, W}^{(I_1, \cdots, I_k)} = \mathcal{R}^* \big( (\boxtimes_j \mathcal{S}_{I_j, W_j}) \boxtimes \mathcal{A} \big)$. The argument in $loc.cit.$ still works.

\subsection{Fusion for nearby cycles on the curve} \label{subsectionFusiononCurve}

To simplify the notation, in this subsection we assume that $J$ is empty. It is easy to generalize the result to general $J$.

\subsubsection{Appplication of Drinfeld's lemma}  \label{subsection-ES-relations-example}
For $\mathcal{F}_{I, W} = \mathcal{R}^* (\mathcal{S}_{I, W} \boxtimes \mathcal{A})$, for any $v \in X$ and $V \in \mathrm{Rep}_{\Lambda} {}^LG$, the construction of excursion operator $S_{V, v}$ in \protect{\cite[Section 6]{Lafforgue}} still works for $\mathfrak{p}_! \mathcal{F}_{I, W}$. We still have $T(h_{V, v}) = S_{V, v}$ for $v \in X-x$. We still have the Eichler-Shimura relations (\protect{\cite[Section 7]{Lafforgue}} for $\Lambda = E$, \protect{\cite[Section 7.2]{XueIntegral}} for $\Lambda = \mathcal{O}_E, k_E$) for $\mathfrak{p}_! \mathcal{F}_{I, W}$.

\begin{prop}   \label{prop-p-F-I-V-constant}
(1) For all $(I,W)$, the complex $\mathfrak{p}_!\mcF_{I,W}$ is constant over over $(\overline{\eta})^{I}$. 

(2) Its geometric generic fiber is equipped with an action of $\mathrm{Weil}(\eta, \overline{\eta})^{I}$. 
\end{prop}
\begin{proof}
This is a mild generalization of \cite{XueSmoothness}.
By the above discussion the complex $\mathfrak{p}_!\mcF_{I,W}$ is equipped with an action of the partial Frobenius morphisms and has the Eichler-Shimura relations. So the argument in \cite{XueSmoothness} to apply Drinfeld's lemma to each $R^j \mathfrak{p}_!\mcF_{I,W}$, $j 
\in \mathbb{Z}$ still works. We deduce that $R^j \mathfrak{p}_!\mcF_{I,W}$ is a constant sheaf over $(\overline{\eta})^{I}$. 
We conclude by Lemma \ref{lem-complexe-is-constant} below (applied to $Y=(\overline{\eta})^{I}$). 
\end{proof}

\begin{lem}  \label{lem-complexe-is-constant}
Let $Y$ be a scheme over $\Fqb$ which does not have cohomology (i.e. denote by $\pi: Y \rightarrow \mathrm{Spec}(\Fqb)$, then the canonical unit map $\mathrm{Id} \rightarrow \pi_* \pi^*$ is an isomorphism). Let $\mathcal{K} \in D_c^{(-)}(Y)$. If $\mathcal{K}^j$ is constant for every $j \in \mathbb{Z}$, then $\mathcal{K}$ is constant.
\end{lem}
\begin{proof}
It is obvious that the following statements are equivalent:
\begin{enumerate}
\item the complex $\mathcal{K}$ is constant 
\item the canonical co-unit map $\pi^* \pi_* \mathcal{K} \rightarrow \mathcal{K}$ is an isomorphism 
\item for every $i \in \mathbb{Z}$, the induced degree $i$ map $R^i \pi^* \pi_* \mathcal{K} \rightarrow \mathcal{K}^i$ is an isomorphism. 
\end{enumerate}

We prove the statement (3): since $Y$ does not have cohomology, for every $i \in \mathbb{Z}$, $R^i \pi^* \pi_* \mathcal{K} = \pi^* \pi_*  \mathcal{K}^i$. Since the sheaf $\mathcal{K}^i$ is constant, the canonical map $\pi^* \pi_* \mathcal{K}^i \rightarrow \mathcal{K}^i$ is an isomorphism. So $R^i \pi^* \pi_* \mathcal{K} \rightarrow \mathcal{K}^i$ is an isomorphism.
\end{proof}

\begin{prop}  \label{prop-p-Psi-F-I-V-constant}
(1) Let $I = I_1 \cup I_2$ and $W=W_1 \boxtimes W_2$. The complex $\mathfrak{p}_! \Psi_{\overline{\eta}_{I_2} \to x} \mcF_{I,W}$ is constant over $(\overline{\eta})^{I_1}$.

(2) Its geometric generic fiber is equipped with an action of $\mathrm{Weil}(\eta, \overline{\eta})^{I_1}$.
\end{prop}

\begin{proof}
By Lemma \ref{lem-Psi-F-equals-F-prime} below, 
$\mathfrak{p}_! \Psi_{\overline{\eta}_{I_2} \to x} \mcF_{I,W} = \mathfrak{p}_! \mcF_{I_1,W_1}'$. We deduce the result from Proposition \ref{prop-p-F-I-V-constant}. 
\end{proof}

\begin{rque}
Another way to prove Proposition \ref{prop-p-Psi-F-I-V-constant} (2) which does not use Lemma \ref{lem-Psi-F-equals-F-prime}: since partial Frobenius morphisms are homeomorphisms, they commute with nearby cycles. We deduce that $\mathfrak{p}_! \Psi_{\overline{\eta}_{I_2} \to x} \mcF_{I,W}$ is equipped with an action of the partial Frobenius morphisms. Then we prove that the Eichler-Shimura relations still hold for $\mathfrak{p}_! \Psi_{\overline{\eta}_{I_2} \to x} \mcF_{I,W}$ and apply Drinfeld's lemma.
\end{rque}

\begin{lem}  \label{lem-Psi-F-equals-F-prime}
Let $I = I_1 \cup I_2$ and $W=W_1 \boxtimes W_2$. Let $\mathcal{F}_{I, W} = \mathcal{R}^*(\mathcal{S}_{I, W} \boxtimes \mathcal{A})$ with $\mathcal{A} \in D_c^b(\ChtR^{nx}_{\{0\}, W_0} |_{x}, \Lambda)$. Then 
$$\Psi_{\overline{\eta}_{I_2} \to x} \mcF_{I,W} = \mathcal{F}'_{I_1, W_1} $$
where $\mathcal{F}'_{I_1, W_1}$ is another sheaf of the form $\mathcal{R}^*(\mathcal{S}_{I_1, W_1} \boxtimes \mathcal{B})$, for some $W_0' \in \Rep_{\Lambda} {^L}G$ and sheaf $\mathcal{B} \in D_c^b(\ChtR^{nx}_{\{0\}, W_0'} |_{x}, \Lambda)$. 
\end{lem}

\begin{proof}
%The nearby cycle $\Psi_{\overline{\eta}_{I_2} \to x}$ is taken with respect to the projection to $X^{I_2}$.
Consider stacks of shtukas with intermediate modifications and the smooth morphism $$\Cht_{I_1 \cup I_2 \cup \{0\}, W_1 \boxtimes W_2 \boxtimes W_0}^{(I_1, I_2, 0)} \xrightarrow{\mathcal{R}} \ChtR_{I_1 \cup I_2 \cup \{0\}, W_1 \boxtimes W_2 \boxtimes W_0}^{(I_1, I_2, 0)}$$
%Note that $\mathcal{R}$ is smooth, so $\mathcal{R}^*$ commutes with the nearby cycle functor. 
By Lemma \ref{lemTwistedProductStructure}, on the restriction over $(\overline{\eta})^{I_1} \times \overline{\eta}_{I_2} \times x$, we have $\mathcal{F}_{I, W}^{(I_1, I_2)} = \mathcal{R}^*(\mathcal{S}_{I_1, W_1}  \widetilde{\boxtimes} \mathcal{S}_{I_2, W_2} \boxtimes \mathcal{A} )$.
We deduce
\begin{equation}   \label{equation-Psi-3-F-1-2-3}
\begin{aligned}
\Psi_{\overline{\eta}_{I_2} \to x}  \mathcal{F}_{I, W}^{(I_1, I_2)}  & = \Psi_{\overline{\eta}_{I_2} \to x}  \mathcal{R}^*(\mathcal{S}_{I_1, W_1}   \widetilde{\boxtimes} \mathcal{S}_{I_2, W_2} \boxtimes \mathcal{A} ) \\
& \simeq \mathcal{R}^* \Psi_{\overline{\eta}_{I_2} \to x}  (\mathcal{S}_{I_1, W_1}  \widetilde{\boxtimes} \mathcal{S}_{I_2, W_2} \boxtimes \mathcal{A} ) \\
& \simeq \mathcal{R}^*  (\mathcal{S}_{I_1, W_1}  \boxtimes \Psi_{\overline{\eta}_{I_2} \to x}  (\mathcal{S}_{I_2, W_2} \boxtimes \mathcal{A} ) )
\end{aligned}
\end{equation}

We have the Cartesian diagram
$$
\xymatrixrowsep{1.5pc}
\xymatrixcolsep{5pc}
\xymatrix{
\Cht_{I \cup \{0\}}^{(I_1, I_2, 0)} |_{(\overline{\eta})^{I_1} \times x^{I_2} \times x}  \ar[r]^{\mathcal{R} \quad \quad }  \ar[d]^{\pi}
& [L^+_{I_1}G \backslash \Gr_{I_1}] |_{(\overline{\eta})^{I_1}}   \times \ChtR_{I_2 \cup \{0\}}^{nx, (I_2, 0)}|_{x^{I_2} \times x}  \ar[d]^{\pi}  \\
\Cht_{I \cup \{0\}} |_{(\overline{\eta})^{I_1} \times x^{I_2} \times x}  \ar[r]^{\mathcal{R} \quad \quad}
& [L^+_{I_1}G \backslash \Gr_{I_1}]|_{(\overline{\eta})^{I_1}}   \times \ChtR_{I_2 \cup \{0\}}^{nx} |_{x^{I_2} \times x}
}
$$

Since $\pi$ is proper, $\pi_!$ commutes with $\Psi_{\overline{\eta}_{I_2} \to x}$. Applying $\pi_!$ to (\ref{equation-Psi-3-F-1-2-3}), we deduce that 
$$\Psi_{\overline{\eta}_{I_2} \to x}  \mathcal{F}_{I, W}  = \mathcal{R}^*(\mathcal{S}_{I_1, W_1}  \boxtimes \mathcal{B})$$
over $$\Cht_{I_1 \cup I_2 \cup \{0\}, W_1 \boxtimes W_2 \boxtimes W_0} |_{(\overline{\eta})^{I_1} \times x^{I_2} \times x} = \Cht_{I_1 \cup \{0\}, W_1 \boxtimes (W_2^{\xi} \otimes W_0)} |_{(\overline{\eta})^{I_1} \times x} $$ 
with $\mathcal{B} = \pi_! \Psi_{\overline{\eta}_{I_2} \to x} (\mathcal{S}_{I_2, W_2} \boxtimes \mathcal{A})$ over $\ChtR_{I_2 \cup \{0\}, W_2 \boxtimes W_0}^{nx} |_{x^{I_2} \times x} =\ChtR^{nx}_{\{0\}, W_2^{\xi} \otimes W_0}  |_x$, where $W_2^{\xi}$ is $W_2$ viewed as representation of ${^L}G$ via the diagonal action.

We take $\mathcal{F}'_{I_1, W_1} = \mathcal{R}^*(\mathcal{S}_{I_1, W_1}  \boxtimes \mathcal{B})$.
\end{proof}

\subsubsection{Fusion properties}

For all surjective maps of finite sets $\xi : I \to K$, we denote by $$\Delta_{\xi}: X^K \to X^I, (x_j)_{j \in K} \mapsto (x_i)_{i \in I}, x_i = x_{\xi(i)}$$ the morphism induced by $\xi$.

\begin{lem}\label{lemFusionNearbyCycleCohomSheaves}
Let $\mathcal{G}$ be a constant complex over $(\overline{\eta})^{I}$. Then for all surjective maps of finite sets $\xi : I \to K$ and all partitions $I = I_1 \cup \dots \cup I_k$, there are canonical isomorphisms
\begin{enumerate}
\item $\Psi_{\overline{\eta}_K \to x}\Delta_{\xi}^*\mcG = \Psi_{\overline{\eta}_I \to x} \mcG$,
\item $\Psi_{\overline{\eta}_I \to x} \mcG = \Psi_{\overline{\eta}_{I_1} \to x}\dots \Psi_{\overline{\eta}_{I_k} \to x}\mcG$.
\end{enumerate}
\end{lem}

\begin{proof}
(1) The first isomorphism comes from the base change property along $\Delta_{\xi}: X^K \to X^I$. We have $\Psi_{\overline{\eta}_K \to x}\Delta_{\xi}^*\mcG \simeq \Psi_{\Delta_{\xi}(\overline{\eta}_K) \to x} \mcG$. Since $\mcG$ is constant, we have canonical isomorphism $\Psi_{\Delta_{\xi}(\overline{\eta}_K) \to x} \mcG \simeq \Psi_{\overline{\eta}_I \to x} \mcG$.

(2) $\mathcal{G}$ is of the form $\boxtimes_{I} \mathcal{A}_i$ and the second isomorphism comes from the Künneth formula for nearby cycles. 

(Another way to see is that all the above items are canonically isomorphic to the fiber $\otimes_{I} (\mathcal{A}_i)_{\overline{\eta}}$.)
\end{proof}

In \ref{subsectionConstructionInverse}, we will apply this lemma to $\mathcal{G} = \mathfrak{p}_!\mcF_{I,W}$ or $\mathcal{G} = \mathfrak{p}_! \Psi_{\overline{\eta}_{I_2} \to x} \mcF_{I,W}$.

\subsection{Fusion of nearby cycles on shtukas} \label{subsectionFusiononshtukas}

We denote by $\mcY \to X^I$ a stack locally of finite type and by $j : \mcY_{\check{X}^I} \to \mcY$ the corresponding open inclusion. We also fix $\mcF \in \DD_c^b(\mcY_{\check{X}^I}, \Lambda)$. The reader should keep in mind that we will apply the construction of this section to the stack $\Cht_{I \cup \{0\}}$  and $\mcF = \mcF_{I,W}$ in the later sections.

\subsubsection{Strict transform of full flags.}  \label{subsection-stricttransform}
We first explain a strict transform construction as in \cite[Remarque 3.2]{GenestierLafforgue}. We consider the following setting 
\begin{enumerate}
\item let $S$ be an integral normal scheme of finite type, 
\item let $f : T \to S$ be a modification, 
\item let $D \subset S$ be an integral Cartier divisor. 
\end{enumerate}
In this situation, proceeding essentially as in \cite[Tag 080C]{Stacks} we define the strict transform $\tilde{D}$ of $D$. We denote by $\eta_D$ the generic point of $D$. 

\begin{lem}\label{lemStrictTransform}
There is a unique point in the fiber of $f$ over $\eta_D$. 
\end{lem}
\begin{proof}
Consider the natural map $\Spec(\mcO_{S, \eta_D}) \to S$ and the pullback square
\[\begin{tikzcd}
	{T \times_S \Spec(\mcO_{S, \eta_D})} & T \\
	{\Spec(\mcO_{S, \eta_D})} & S.
	\arrow[from=1-1, to=1-2]
	\arrow[from=1-1, to=2-1]
	\arrow[from=1-2, to=2-2]
	\arrow[from=2-1, to=2-2]
\end{tikzcd}\]
Since $S$ is normal and $D$ has codimension one $\mcO_{S, \eta_D}$ is DVR. The map $T \times_S \Spec(\mcO_{S, \eta_D}) \to \Spec(\mcO_{S, \eta_D})$ is still proper and surjective. The generic point of both the source and the target are respectively the generic points of $S$ and $T$ hence this map is also generically an isomorphism and is thus a modification. But there are no nontrivial modifications of a valuative ring and so this map is an isomorphism. The lemma follows.
\end{proof}

We denote by $\eta_{\tilde{D}}$ the unique point above $\eta_D$ and by $\tilde{D} \subset T$ its reduced closure. We call $\tilde{D}$ the strict transform of $D$ along $f$. The proof lemma \ref{lemStrictTransform} also immediatly yields the following corollary.

\begin{corol}
Let $T_1 \to T_2 \to S$ be a tower of modifications and denote by $\tilde{D}_1$ and $\tilde{D}_2$ the strict transform of $D$ along $T_1 \to S$ and $T_2 \to S$. Then $\tilde{D}_1$ is also the strict transform of $\tilde{D}_2$ along $T_1 \to T_2$. 
\end{corol}

\begin{lem}\label{lemModifStrictTransform}
The map $\tilde{D} \to D$ is a modification. 
\end{lem}

\begin{proof}
By construction this map is generically an isomorphism and as both its source and target are irreducible, it is dominant. Since it factors as $\tilde{D} \to f^{-1}(D) \to D$ and the first map is a closed immersion this map is proper. Therefore it is also surjective and thus it is a modification. 
\end{proof}

We consider now the following situation : 
\begin{enumerate}
\item $S$ is as before an integral normal scheme of finite type, 
\item $f : T \to S$ is a modification, 
\item $S_0 \subset S_{1} \subset \dots \subset S_n = S$ is sequence of integral normal closed subschemes of $S$ such that $S_{i-1}$ is a Cartier divisor of $S_{i}$. 
\end{enumerate}
In this situation, we define $T_{i-1}$ inductively to be the strict transform of $S_{i-1}$ along $T_{i} \to S_{i}$. In the case where $S_0$ is a point, then so is $T_0$. As in \cite{GenestierLafforgue}, we call the point $T_0$ the strict transform of $S_0$ along the flag $(S_i)$. 

Finally we consider the situation where $S = X^n$ and we define the following full flags of $S$.
\begin{enumerate}
\item The diagonal flag : 
\begin{equation}\label{diagonalFlag}
x \in X \subset X^2 \subset \dots \subset X^{n-1} \subset X^n,
\end{equation}
where each incusion $X^i \subset X^{i+1}$ is the partial diagonal along the last two coordinates, namely this is the map 
$$(x_1, \dots, x_{i-1}) \mapsto (x_1, \dots, x_{i-1}, x_{i-1}).$$
\item Let $d \in \N^n$ be a multiset and recall that we denote by $\Frob^d$ the corresponding partial Frobenius morphism. We denote by 
$$x = Y^d_0 \subset Y_1^d = \Frob^d(X) \subset Y^d_2 = \Frob^d(X^2) \dots \subset Y_n^d = X^n$$
the image of the diagonal under the partial Frobenius morphism $\Frob^d$.
\item The hyperplane flag : 
\begin{equation}\label{hyperplaneFlag}
x = x^n \in X \times x^{n-1} \subset \dots \subset X^{n-1} \times x \subset X^n,
\end{equation}
this is the flag for which the $k$-th term is $X^{k}$ and inclusion in $X^n$ sets the last $n-k$ coordinates to $x$. 
\end{enumerate}

The next lemma is \cite[Lemme 3.3]{GenestierLafforgue}.
\begin{lem}\label{lemStrictTransformFrob}
Let $\widetilde{S} \to X^n$ be a modification, there exists $d \in \N^n$ large enough and increasing enough such that the strict transform of $x$ along the hyperplane flag and the strict transform of $x$ along the flag $(Y_i^d)$ are the same point $s \in \widetilde{S}$.
\end{lem}

Let $I$ be a finite set, given a total ordering of $I$, there is a unique increasing bijection $I = \{1, \dots, n\}$. In the rest of this section, we will work with general totally ordered sets $I$, the construction of the two previous flags then naturally extends to this setting. 

\subsubsection{Restricting to partial diagonal}\label{subsubsectionPartialDiagonal}

Let $f : \widetilde{S} \to X^I$ be a modification such that $\mcF$ is $(\Psi, \check{X}^I)$-good (we refer to Definition \ref{def-Psigood}).
We denote by $(\widetilde{S}_{i, 0})$ the strict transform of the diagonal flag (\ref{diagonalFlag}) in $\widetilde{S}$ and by $s_0$ the strict transform of $x$. More generally, let $d \in \N^I$, we denote by $(\widetilde{S}_{i,d})$ the strict transform of the flag $(Y_i^d)$ in $\widetilde{S}$ and by $s_d$ the strict transfrom of $x$ along this flag. 

We denote by $\eta_{i, d}$ the generic point of $\widetilde{S}_{i, d}$ and we fix a chain of specializations of geometric points
\begin{equation}  \label{equation-specializations-Frob-diag}
\overline{\eta}_I = \overline{\eta}_{n, d} \to \overline{\eta}_{n-1, d} \to \dots \to s_d.
\end{equation}

\begin{lem}   \label{lemGabberThm}
Suppose that $\mcF$ is ULA relative to $\check{X}^I$. 
\begin{enumerate}
\item For any $d$, for all $1 \leq i \leq n$, the natural map 
$$\Psi_{\overline{\eta}_I \to s_d}^{\widetilde{S}} \mcF \to \Psi_{\overline{\eta}_{i, d} \to s_d}^{\widetilde{S}_{i, d}}\mcF$$
is an isomorphism. 
\item For any $d$, for all $1 \leq i' \leq i \leq n$, the natural following triangle is commutative and all maps are isomorphisms.
\[\begin{tikzcd}
	{\Psi_{\overline{\eta}_I \to s_d}^{\widetilde{S}}\mcF} \\
	{\Psi_{\overline{\eta}_{i, d} \to s_d}^{\widetilde{S}_{i, d}}\mcF} & {\Psi_{\overline{\eta}_{i', d} \to s_d}^{\widetilde{S}_{i', d}}\mcF}
	\arrow[from=1-1, to=2-1]
	\arrow[from=1-1, to=2-2]
	\arrow[from=2-1, to=2-2]
\end{tikzcd}\] 
\end{enumerate}
\end{lem} 

\begin{proof}
$(1)$. The first map is the following composition 
\begin{align*}
\Psi_{\overline{\eta}_I \to s_d}^{\widetilde{S}}\mcF &\to \Psi_{\overline{\eta}_{i, d} \to s_d}^{\widetilde{S}}\Psi_{\overline{\eta}_I \to \overline{\eta}_{i, d}}^{\widetilde{S}} \mcF \\
&\simeq \Psi_{\overline{\eta}_{i, d} \to s_d}^{\widetilde{S}} \mcF \\
&\to \Psi_{\overline{\eta}_{i, d} \to s_d}^{\widetilde{S}_{i, d}}\mcF.
\end{align*}
The first map is the map coming from Gabber's theorem \ref{thmGabber} and is an isomorphism by the $\Psi$-goodness hypothesis. The second map follows from the observation that the two points $\overline{\eta}_I$ and $\overline{\eta}_{i, d}$ lie in $f^{-1}(\check{X}^I)$, the open subset of $\widetilde{S}$ where $\mcF$ is ULA. For ULA sheaves the nearby cycle functor is isomorphic to the identity functor. Finally, the last map is an isomorphism by Orgogozo's theorem \ref{thmOrgogozo}. The point $(2)$ is an iteration of the same argument. 
\end{proof}

\begin{rque}
Both maps $\Psi_{\overline{\eta}_I \to s_d}\mcF \to \Psi_{\overline{\eta}_{i, d} \to s_d}\mcF$ and $\Psi_{\overline{\eta}_{i ,d} \to s_d}\mcF \to \Psi_{\overline{\eta}_{i', d} \to s_d}\mcF$ are constructed using the base change maps of higher nearby cycles and the canonical maps for the compositions of nearby cycles hence their formation commute with pullback along smooth maps, pushforward along proper maps and are in general equipped with canonical base change maps against $!$-pushforward.
\end{rque}

\subsubsection{Iterated nearby cycles}\label{subsubsectionIteration}

We fix a total order on $I = \{1, \dots, n\}$ and we will equip $X^I$ with the hyperplane flag (\ref{hyperplaneFlag}). Now we construct a modification $\widetilde{S}$ of $X^I$:

Let $\xi : I \to K = \{1, \dots, k\}$ be a surjective map of totally ordered finite sets, we denote by $I_j = \xi^{-1}(j)$ and by $I_{\geq j} = \xi^{-1}(\{j' \geq j\})$. We construct by induction two modifications $\widetilde{T}_j$ and $T_j$ of $X^{I_j}$ and $X^{I_{\geq j}}$. 
\begin{enumerate}
\item We let $\widetilde{T}_k = T_k$ be a modification of $X^{I_k}$ such that $\mcF$ is $(\Psi, \check{X}^{I_k})$-good relative to $\mathcal{Y} \to X^I \to X^{I_k}$ where the second map is induced by the projection along $I_k \subset I$. 
\item Let $j \in K$ and assume that for all $j' > j$ the modifications $T_{j'}$ and $\widetilde{T}_{j'}$ have been constructed. The total orders on $I_{j'}$ determine points $t_{j'} \in \widetilde{T}_{j'}$ over $x$. We let $\widetilde{T}_j \to X^{I_j}$ be a modification such that $\Psi_{\overline{\eta_{I_{j+1}}} \to t_{j+1}}\dots \Psi_{\overline{\eta_{I_k}} \to t_k}\mcF$ is $\Psi$-good relative to the projection to $X^{I_j}$. And we let $T_j \to \widetilde{T}_j \times T_{j+1}$ be a modification such that $\mcF$ is $\Psi$-good relative to the projection to $X^{I_{\geq j}}$. 
\end{enumerate}
Finally we denote by $\widetilde{S} = T_1$. 

We denote by $\widetilde{S}_i$ the strict transform of the $i$-th term of the hyperplane flag (\ref{hyperplaneFlag}) and by $s$ the strict transform of $x$ (which is the $0$-th term of the flag).

Denote by $i_j \in I$ the maximal element in $I_j$. We denote by $\eta_i$ the generic point of $\widetilde{S}_i$ the $i$-th term in the flag $(\widetilde{S}_i)$. So that $\eta_{|I|} = \eta_{i_k} = \eta_I$. We also fix specialization maps
\begin{equation} \label{equation-specializations-hyperplan}
\overline{\eta}_n \to \overline{\eta}_{n-1} \dots \to s
\end{equation}  
By Gabber's and Orgogozo's theorem \ref{thmGabber} and \ref{thmOrgogozo}, we have a canonical isomorphism 
\begin{equation}  \label{equation-Gabber-Orgogozo-hyperplan}
\begin{aligned}
\Psi_{\overline{\eta}_I \to s}^{\widetilde{S}}\mcF &\simeq \Psi_{\overline{\eta}_{i_1} \to s}^{\widetilde{S}} \Psi_{\overline{\eta}_{i_2} \to \overline{\eta}_{i_1}}^{\widetilde{S}} \dots \Psi_{\overline{\eta}_{i_k} \to \overline{\eta}_{i_{k-1}}}^{\widetilde{S}}\mcF \\
&\simeq \Psi_{\overline{\eta}_{i_1} \to s}^{{\widetilde{S}}_{i_1}} \Psi_{\overline{\eta}_{i_2} \to \overline{\eta}_{i_1}}^{{\widetilde{S}}_{i_2}} \dots \Psi_{\overline{\eta}_{i_k} \to \overline{\eta}_{i_{k-1}}}^{\widetilde{S}_{i_k}}\mcF
\end{aligned}
\end{equation}
Finally consider the map $\widetilde{S}_{i_j} \to \widetilde{T}_{j}$, this map is surjective and sends the specialization map $\overline{\eta}_{i_j} \to \overline{\eta}_{i_{j-1}}$ to the specialization map $\overline{\eta}_{I_j} \to t_j$ hence we have a canonical base change map (Remark \ref{rquebasechange}) (which need not be an isomorphism):
\begin{equation}\label{equationKunnethMaps}
\Psi_{\overline{\eta}_{I_1} \to t_1}^{\widetilde{T}_1} \Psi_{\overline{\eta}_{I_2} \to t_2}^{\widetilde{T}_2} \dots \Psi_{\overline{\eta}_{I_k} \to t_k}^{\widetilde{T}_k}\mcF \to \Psi_{\overline{\eta}_{i_1} \to s}^{\widetilde{S}_{i_1}} \Psi_{\overline{\eta}_{i_2} \to \overline{\eta}_{i_1}}^{\widetilde{S}_{i_2}} \dots \Psi_{\overline{\eta}_{i_k} \to \overline{\eta}_{i_{k-1}}}^{\widetilde{S}_{i_k}}\mcF.
\end{equation}
Both maps are compatible with pushforward along proper maps and pullback along smooth maps and are equipped with base change maps against $!$-pushforward.

\begin{lem}\label{lemIteratedNearbyCyclesOnshtukas}
Let $\mcF = \mcF_{I,W}$ and $\mcY = \Cht_{I \cup J \cup \{0\}}$. We fix a map $I \to K= \{1, \dots, k\}$ of totally ordered sets and we assume that $W = W_1 \boxtimes \dots \boxtimes W_k$.  
Then the map (\ref{equationKunnethMaps}) is an isomorphism.  
\end{lem}

\begin{proof}
Consider the partition $I = I_1 \sqcup \dots \sqcup I_k$ induced by the map $I \to K$.  
Recall (see \ref{subsection-ChtR-intermediate}) that we have a proper map $\pi: \ChtR^{nx, (I_1, \dots, I_k, J \cup \{0\})}_{I \cup J \cup \{0\}, (n_i)} \to \ChtR_{I \cup J \cup \{0\}, (n_i)}^{nx}$ and $\mcF_{I,W}= \pi_! \mcF_{I,W}^{(I_1, \dots, I_k)}$. Since $\pi$ is proper, $ \pi_!$
commutes with $\Psi$, hence it is enough to prove that the map (\ref{equationKunnethMaps}) is an isomorphism upstairs. By lemma \ref{lemTwistedProductStructure}, the sheaf upstairs $\mcF_{I,W}^{(I_1, \dots, I_k)}$ is a twisted product.
Arguing as in lemma \ref{lem-two-legs-Psi-delta-d-eta-equal-Psi-1-Psi-2}, step (4), by the Künneth formula for nearby cycles \ref{thmKunnethCyclesProches} we get 
$$ \Psi_{\overline{\eta}_{I_k} \to t_k}^{\widetilde{T}_k}\mcF_{I,W}^{(I_1, \dots, I_k)} \xrightarrow{\sim} \Psi_{\overline{\eta}_{i_k} \to \overline{\eta}_{i_{k-1}}}^{\widetilde{S}_{i_k}}\mcF_{I,W}^{(I_1, \dots, I_k)}.$$
Again both sides are twisted products. 
Applying the Künneth formula successively, we get the desired isomorphism. 
\end{proof}

\subsubsection{Fusion for nearby cycles on shtukas}

Firstly we introduce two copies of $I$ which we denote by $I_1$ and $I_2$. Let $W_1, W_2 \in \Rep {^L}G^I$ and consider the sheaf $\mcF = \mcF_{I_1 \cup I_2, W_1 \boxtimes W_2}$ which is a sheaf on $\Cht_{I_1 \cup I_2 \cup J \cup \{0\}}$ over $(\check{X})^{I_1 \cup I_2 \cup J} \times x$. Let $\Delta^{1,2} : X^I \to X^{I_1 \cup I_2}$ be the diagonal map. For any $d = (d_1, d_2) \in \N^{I_1 \cup I_2}$, we denote by $\Delta^{1, 2}_d := F^{d}(\Delta^{1, 2})$.

\begin{lem}\label{lemFusionOnshtukas}
There exists a modification $\widetilde{S} \to X^{I_1} = X^{I_2}$, a point $s \in \widetilde{S}$, a modification $\widetilde{S}_{12} \to \widetilde{S} \times \widetilde{S}$, a point $s_{12} \in \widetilde{S}_{12}$ above $(s, s)$ and a tuple of integer $d \in \N^{I_1 \cup I_2}$ such that there is a canonical isomorphism
$$\Psi_{\Delta^{1,2}_d(\overline{\eta}_{I}) \to s_{12}} \mcF = \Psi_{\overline{\eta}_{I_1 \cup I_2} \to s_{12}}\mcF =  \Psi_{\overline{\eta}_{I_1} \to s} \Psi_{\overline{\eta}_{I_2} \to s}\mcF.$$
\end{lem}

\begin{proof}
We fix a total order on $I$ which determines an order on $I_1 \cup I_2$ by declaring that for all $i \in I_1$ and $j \in I_2$ we have $i < j$.
Applying the construction \ref{subsubsectionIteration} to the partition given by $(I_1, I_2)$, we get two modifications $\widetilde{X}^{I_i} \to X^{I_i}$ and a modification $\widetilde{S}_{12} \to \widetilde{X}^{I_1} \times \widetilde{X}^{I_2}$, when choosing these modifications we can further assume that $\widetilde{X}^{I_1}$ dominates $\widetilde{X}^{I_2}$ and that $\widetilde{S}_{12}$ dominates $\widetilde{X}^{I_1} \times \widetilde{X}^{I_2}$. Let $\widetilde{S} = \widetilde{X}^{I_1}$.
The points $s \in \widetilde{S}$ and $s_{12} \in \widetilde{S}_{12}$ are the strict transforms of $x$ along the hyperplane flag determined by the order on $I_1$ and $I_1 \cup I_2$. 
Let $\eta_{I_1}$ be the generic point of the strict transform of $X^{I_1} \times x$ in $\widetilde{S}_{12}$.
By construction \ref{subsubsectionIteration} we have isomorphism
\begin{equation}
\Psi_{\overline{\eta}_{I_1 \cup I_2} \to s_{12}}\mcF \simeq \Psi_{\overline{\eta}_{I_1} \to s_{12}} \Psi_{\overline{\eta}_{I_1 \cup I_2} \to \overline{\eta}_{I_1}}\mcF \simeq \Psi_{\overline{\eta}_{I_1} \to s} \Psi_{\overline{\eta}_{I_2} \to s}\mcF
\end{equation}
where the first isomorphism is (\ref{equation-Gabber-Orgogozo-hyperplan}) and the second isomorphism is (\ref{equationKunnethMaps}) and lemma \ref{lemIteratedNearbyCyclesOnshtukas}.

By lemma \ref{lemStrictTransformFrob}, we can choose $d \in \N^{I_1 \cup I_2}$ increasing enough such that the strict transform $(s_{12})_d \in \widetilde{S}_{12}$ of $x$ along the flag $(Y_i^d)$ and 
the above strict transform $s_{12} \in \widetilde{S}_{12}$ of $x$ along the hyperplane flag agree. Choose specialization maps in (\ref{equation-specializations-Frob-diag}) such that the composition coincides with the composition of (\ref{equation-specializations-hyperplan}). 
By construction \ref{subsubsectionPartialDiagonal} Lemma \ref{lemGabberThm}, we have isomorphism
\begin{equation*}
\Psi_{\Delta^{1,2}_d(\overline{\eta}_{I}) \to s_{12}} \mcF \simeq \Psi_{\overline{\eta}_{I_1 \cup I_2} \to s_{12}}\mcF 
\end{equation*}
\end{proof}

\begin{rque}
The only place where we need $d$ large enough is to relate the strict transform of $x$ along the flag $(Y_i^d)$ and the strict transform of $x$ along the hyperplane flag, thus relate construction \ref{subsubsectionPartialDiagonal} and construction \ref{subsubsectionIteration}. 
\end{rque}

The next lemma is proved as in Lemma \ref{lem-p-Frob-d-Delta-equal-p-Delta}.
\begin{lem}\label{lemPartialFrobResDiag}
We keep the notations as in lemma \ref{lemFusionOnshtukas}. For any $d \in \N^{I_1 \cup I_2}$, the partial Frobenius map induces an isomorphism
\begin{equation*}  
\mathfrak{p}_! \Psi_{\overline{\eta}_I \to s}\Delta^{1,2,*}_d \mcF = \mathfrak{p}_!\Psi_{\overline{\eta}_I \to s}\Delta^{1,2,*}\mcF.
\end{equation*}
\end{lem}

\begin{lem}\label{lemFusionNearbyCycle2}
There exists a modification $\widetilde{S} \to X^{I_1} = X^{I_2}$ and a point $s \in \widetilde{S}$ such that there is a canonical isomorphism 
\begin{equation*}
\mathfrak{p}_!\Psi_{\overline{\eta}_I \to s}\Delta^{1,2,*}\mcF \to \mathfrak{p}_!\Psi_{\overline{\eta}_{I_1} \to s}\Psi_{\overline{\eta}_{I_2} \to s}\mcF.
\end{equation*}
\end{lem}

\begin{proof}
We have the following isomorphisms:
$$\mathfrak{p}_!\Psi_{\overline{\eta}_I \to s}\Delta^{1,2,*}\mcF \simeq \mathfrak{p}_! \Psi_{\overline{\eta}_I \to s}\Delta^{1,2,*}_d \mcF \simeq \mathfrak{p}_!\Psi_{\Delta_{d}^{1,2}(\overline{\eta}_I) \to s}\mcF \simeq \mathfrak{p}_!\Psi_{\overline{\eta}_{I_1} \to s}\Psi_{\overline{\eta}_{I_2} \to s}\mcF$$
The first isomorphism is by Lemma \ref{lemPartialFrobResDiag}. The second isomorphism is by the base change (Remark \ref{rquebasechange}). The last isomorphism is by lemma \ref{lemFusionOnshtukas}.
\end{proof}

Now we introduce three copies of $I$ which we denote by $I_1, I_2$ and $I_3$. Let $W \in \Rep {^L}G^I$ and consider the sheaf 
$\mcG = \mcF_{I_1 \cup I_2 \cup I_3, W \boxtimes W^* \boxtimes W}$ on $\Cht_{I_1 \cup I_2 \cup I_3 \cup J \cup \{0\}}$ over $(\check{X})^{I_1 \cup I_2 \cup I_3 \cup J} \times x$. We will denote some partial diagonals 
\begin{enumerate}
\item $\Delta^{1,2} : X^I \to X^{I_1 \cup I_2}$, 
\item $\Delta^{1,2} : X^{I_1 \cup I_3} \to X^{I_1 \cup I_2 \cup I_3}$, 
\item $\Delta^{2,3} : X^I \to X^{I_2 \cup I_3}$, 
\item $\Delta^{2,3} : X^{I_1 \cup I_2} \to X^{I_1 \cup I_2 \cup I_3}$, 
\item $\Delta^{1,2,3} : X^I \to X^{I_1 \cup I_2 \cup I_3}$,
\end{enumerate}
obtained in the obvious way. We use the notations $\Delta^{1,2}$ and $\Delta^{2,3}$ to denote two different maps but they should be distinguished from the context.

\begin{lem}\label{lemFusionNearbyCycles3}
There exists a modification $\widetilde{S} \to X^{I_1} = X^{I_2} = X^{I_3}$ and a point $s \in \widetilde{S}$ above $x$ and a commutative diagram of isomorphisms 
\[\begin{tikzcd}
	{\mathfrak{p}_!\Psi_{\overline{\eta}_{I_1} \to s}\Delta^{1,2,3,*}\mcG} & {\mathfrak{p}_!\Psi_{\overline{\eta}_{I_1} \to s}\Psi_{\overline{\eta}_{I_2} \to s}\Delta^{2,3,*}\mcG} \\
	& {\mathfrak{p}_!\Psi_{\overline{\eta}_{I_1} \to s}\Psi_{\overline{\eta}_{I_2} \to s}\Psi_{\overline{\eta}_{I_3} \to s}\mcG} \\
	{\mathfrak{p}_!\Psi_{\overline{\eta}_{I_1} \to s}\Delta^{1,2,3,*}\mcG} & {\mathfrak{p}_!\Psi_{\overline{\eta}_{I_2} \to s}\Psi_{\overline{\eta}_{I_3} \to s}\Delta^{1,2,*}\mcG}
	\arrow[from=1-1, to=1-2]
	\arrow[from=1-1, to=2-2]
	\arrow[Rightarrow, no head, from=1-1, to=3-1]
	\arrow[from=1-2, to=2-2]
	\arrow[from=3-1, to=2-2]
	\arrow[from=3-1, to=3-2]
	\arrow[from=3-2, to=2-2]
\end{tikzcd}\]
making the following diagram commutative 
\[\begin{tikzcd}
	{\mathfrak{p}_!\Psi_{\overline{\eta}_{I_1} \to s}\Delta^{1,2,3,*}\mcG} & {\mathfrak{p}_!\Psi_{\overline{\eta}_{I_1} \to s}\Psi_{\overline{\eta}_{I_2} \to s}\Delta^{2,3,*}\mcG} & {\Psi_{\overline{\eta}_{I_1} \to s}\mathfrak{p}_!\Psi_{\overline{\eta}_{I_2} \to s}\Delta^{2,3,*}\mcG} \\
	& {\mathfrak{p}_!\Psi_{\overline{\eta}_{I_1} \to s}\Psi_{\overline{\eta}_{I_2} \to s}\Psi_{\overline{\eta}_{I_3} \to s}\mcG} & {\Psi_{\overline{\eta}_{I_1} \to s}\mathfrak{p}_!\Psi_{\overline{\eta}_{I_2} \to s}\Psi_{\overline{\eta}_{I_3} \to s}\mcG} \\
	\\
	{\mathfrak{p}_!\Psi_{\overline{\eta}_{I_1} \to s}\Delta^{1,2,3,*}\mcG} & {\mathfrak{p}_!\Psi_{\overline{\eta}_{I_2} \to s}\Psi_{\overline{\eta}_{I_3} \to s}\Delta^{1,2,*}\mcG} & {\Psi_{\overline{\eta}_{I_1} \to s}\Psi_{\overline{\eta}_{I_2} \to s}\mathfrak{p}_!\Psi_{\overline{\eta}_{I_3} \to s}\mcG}
	\arrow[from=1-1, to=1-2]
	\arrow[from=1-1, to=2-2]
	\arrow[Rightarrow, no head, from=1-1, to=4-1]
	\arrow["\can"{description}, from=1-2, to=1-3]
	\arrow[from=1-2, to=2-2]
	\arrow[from=1-3, to=2-3]
	\arrow["\can"{description}, from=2-2, to=2-3]
	\arrow["\can"{description}, from=2-3, to=4-3]
	\arrow[from=4-1, to=2-2]
	\arrow[from=4-1, to=4-2]
	\arrow[from=2-2, to=4-2]
	\arrow["\can"{description}, from=4-2, to=4-3]
\end{tikzcd}\]
where $\can$ is the canonical base change map and the first right vertical map is the map of lemma \ref{lemFusionNearbyCycle2}. 
\end{lem}

\begin{proof}
Denote by $I_{123} = I_1 \cup I_2 \cup I_3$ and by $I_{ij} = I_i \cup I_j$. We apply the construction of modification \ref{subsubsectionIteration} to all the partitions 
\begin{enumerate}
\item $I_{123} = I_{1} \cup I_2 \cup I_3$,
\item $I_{123} = I_{12} \cup I_3$,
\item $I_{123} = I_1 \cup I_{23}$, 
\item $I_{12} = I_{1} \cup I_2$,
\item $I_{23} = I_2 \cup I_3$.
\end{enumerate}
We then get modifications $\widetilde{X}^{I_i}$ of $X^{I_i}$, $\widetilde{S}_{ij}$ of $X^{I_{ij}}$ and $\widetilde{S}_{123}$ of $X^{I_{123}}$. We can assume they have the following dominance relations 
\begin{enumerate}
\item $\widetilde{X}^{I_1}$ dominates $\widetilde{X}^{I_2}$ and $\widetilde{X}^{I_3}$, 
\item $\widetilde{S}_{12} \to \widetilde{X}^{I_1} \times \widetilde{X}^{I_2}$ and $\widetilde{S}_{23} \to \widetilde{X}^{I_2} \times \widetilde{X}^{I_3}$,
\item and $\widetilde{S}_{123} \to \widetilde{X}^{I_1} \times \widetilde{S}_{23}$ and $\widetilde{S}_{123} \to \widetilde{S}_{12} \times \widetilde{X}^{I_3}$. 
\end{enumerate}
We denote by $s_i$ the strict transform of $x$ in $\widetilde{X}^{I_i}$ and by $s_{ij}$ the strict transform of $x$ in $\widetilde{S}_{ij}$ and by $s_{123}$ the strict transform of $x$ in $\widetilde{S}_{123}$. We also denote by $\eta_1 \in \widetilde{S}_{123}$ the generic point of the strict transform of $X^{I_1} \times x$ in $\widetilde{S}_{123}$ and by $\eta_2$ the generic point of the strict transform of $X^{I_1 \cup I_2} \times x$ in $\widetilde{S}_{123}$. 
We choose $d = (d_1,d_2,d_3) \in \N^{I_{123}}$ such that all three tuples $(d_1,d_2,d_3), (d_1,d_2)$ and $(d_2,d_3)$ satisfy the hypothesis of lemma \ref{lemStrictTransformFrob}. This implies that $(s_{123})_d = s_{123} \in \widetilde{S}_{123}$ and $(s_{12})_{d_1,d_2} = s_{12} \in \widetilde{S}_{12}$ and $(s_{23})_{d_2,d_3} = s_{23} \in \widetilde{S}_{23}$. 

Take $\widetilde{S} = \widetilde{X}^{I_1}$ and $s = s_1$. We identify $\Psi_{\overline{\eta}_{I_i} \to s_i}$ and $\Psi_{\overline{\eta}_{I_i} \to s}$.

We now have a diagram of isomorphisms of nearby cycles on the special fiber of $\Cht$ coming from the constructions \ref{subsubsectionIteration} and \ref{subsubsectionPartialDiagonal}, the first diagram of the lemma follows from applying $\mathfrak{p}_!$ and the correct partial Frobenii and looking at the three extremities of the next diagram
\[\begin{tikzcd}
	&&& {\Psi_{\overline{\eta}_{I_{1}} \to s_{1}}\Psi_{\Delta^{2,3}_{(d_2,d_3)}(\overline{\eta}_{I_2}) \to s_{23}}\mcG} \\
	&& {\Psi_{\overline{\eta}_1 \to s_{123}}\Psi_{\overline{\eta}_{I_{123}} \to \overline{\eta}_1}\mcG} & {\Psi_{\overline{\eta}_{I_{1}} \to s_{1}}\Psi_{\overline{\eta}_{I_{23}} \to s_{23}}\mcG} \\
	{\Psi_{\Delta^{1,2,3}_{d}(\overline{\eta}_{I_1}) \to s_{123}}\mcG} & {\Psi_{\overline{\eta}_{I_{123}}\to s_{123}}\mcG} & {\Psi_{\overline{\eta}_1 \to s_{123}}\Psi_{\overline{\eta}_2 \to \overline{\eta}_1}\Psi_{\overline{\eta}_{I_{123}} \to \overline{\eta}_2}\mcG} & {\Psi_{\overline{\eta}_{I_1} \to s_1}\Psi_{\overline{\eta}_{I_2} \to s_2}\Psi_{\overline{\eta}_{I_3} \to s_3}\mcG} \\
	&& {\Psi_{\overline{\eta}_2 \to s_{123}}\Psi_{\overline{\eta}_{I_{123}} \to \overline{\eta}_2}\mcG} & {\Psi_{\overline{\eta}_{I_{12}} \to s_{12}}\Psi_{\overline{\eta}_{I_3} \to s_3}\mcG} \\
	&&& {\Psi_{\Delta^{1,2}_{(d_1,d_2)}(\overline{\eta}_{I_1}) \to s_{12}}\Psi_{\overline{\eta}_{I_3} \to s_3}\mcG}
	\arrow[from=1-4, to=2-4]
	\arrow[from=2-3, to=2-4]
	\arrow[from=2-3, to=3-3]
	\arrow[from=2-4, to=3-4]
	\arrow[from=3-1, to=3-2]
	\arrow[from=3-2, to=2-3]
	\arrow[from=3-2, to=3-3]
	\arrow[from=3-2, to=4-3]
	\arrow[from=3-3, to=3-4]
	\arrow[from=4-3, to=3-3]
	\arrow[from=4-3, to=4-4]
	\arrow[from=4-4, to=3-4]
	\arrow[from=5-4, to=4-4]
\end{tikzcd}\]

For the compatibility with the canonical maps there are two diagrams : 
\[\begin{tikzcd}
	{\mathfrak{p}_!\Psi_{\overline{\eta}_{I_{1}} \to s_{1}}\Psi_{\Delta^{2,3}_{(d_2,d_3)}(\overline{\eta}_{I_2}) \to s_{23}}\mcG} & {\Psi_{\overline{\eta}_{I_{1}} \to s_{1}}\mathfrak{p}_!\Psi_{\Delta^{2,3}_{(d_2,d_3)}(\overline{\eta}_{I_2}) \to s_{23}}\mcG} \\
	{\mathfrak{p}_!\Psi_{\overline{\eta}_{I_{1}} \to s_{1}}\Psi_{\overline{\eta}_{I_{23}} \to s_{23}}\mcG} & {\Psi_{\overline{\eta}_{I_{1}} \to s_{1}}\mathfrak{p}_!\Psi_{\overline{\eta}_{I_{23}} \to s_{23}}\mcG} \\
	{\mathfrak{p}_!\Psi_{\overline{\eta}_{I_1} \to s_1}\Psi_{\overline{\eta}_{I_2} \to s_2}\Psi_{\overline{\eta}_{I_3} \to s_3}\mcG} & {\Psi_{\overline{\eta}_{I_1} \to s_1}\mathfrak{p}_!\Psi_{\overline{\eta}_{I_2} \to s_2}\Psi_{\overline{\eta}_{I_3} \to s_3}\mcG}
	\arrow["\can", from=1-1, to=1-2]
	\arrow[from=1-1, to=2-1]
	\arrow[from=1-2, to=2-2]
	\arrow["\can"{description}, from=2-1, to=2-2]
	\arrow[from=2-1, to=3-1]
	\arrow[from=2-2, to=3-2]
	\arrow["\can"', from=3-1, to=3-2]
\end{tikzcd}\]
and 
\[\begin{tikzcd}
	{\mathfrak{p}_!\Psi_{\overline{\eta}_{I_1} \to s_1}\Psi_{\overline{\eta}_{I_2} \to s_2}\Psi_{\overline{\eta}_{I_3} \to s_3}\mcG} & {\Psi_{\overline{\eta}_{I_1} \to s_1}\mathfrak{p}_!\Psi_{\overline{\eta}_{I_2} \to s_2}\Psi_{\overline{\eta}_{I_3} \to s_3}\mcG} \\
	{\mathfrak{p}_!\Psi_{\overline{\eta}_{I_{12}} \to s_{12}}\Psi_{\overline{\eta}_{I_3} \to s_3}\mcG} & {\Psi_{\overline{\eta}_{I_1} \to s_1}\Psi_{\overline{\eta}_{I_2} \to s_2}\mathfrak{p}_!\Psi_{\overline{\eta}_{I_3} \to s_3}\mcG} \\
	{\mathfrak{p}_!\Psi_{\Delta^{1,2}_{(d_1,d_2)}(\overline{\eta}_{I_1}) \to s_{12}}\Psi_{\overline{\eta}_{I_3} \to s_3}\mcG} & {\Psi_{\Delta^{1,2}_{(d_1,d_2)}(\overline{\eta}_{I_1}) \to s_{12}}\mathfrak{p}_!\Psi_{\overline{\eta}_{I_3} \to s_3}\mcG}
	\arrow["\can", from=1-1, to=1-2]
	\arrow["\can", from=1-2, to=2-2]
	\arrow[from=2-1, to=1-1]
	\arrow["\can"{description}, from=2-1, to=2-2]
	\arrow[from=3-1, to=2-1]
	\arrow["\can"', from=3-1, to=3-2]
	\arrow[from=3-2, to=2-2]
\end{tikzcd}\]
which are commutative by the compatibility of pullbacks for nearby cycles with base change and Künneth maps with base change maps. 
\end{proof}

\subsubsection{The case of a trivial modification}

\begin{lem}\label{lemTrivialModification}
Let $I_1,I_2$ be two finite sets and $W \in \Rep_{\Lambda} {^L}G^{I_1}$, then the sheaf $\mcF_{I_1 \cup I_2, W \boxtimes 1}$ is supported on $\Cht_{I_1 \cup I_2 \cup J \cup \{0\}, W \boxtimes 1 \boxtimes V \boxtimes W_0} = \Cht_{I_1 \cup J \cup \{0\}, W \boxtimes V \boxtimes W_0} \times (\check{X})^{I_2}$ and is isomorphic to $\mcF_{I_1, W} \boxtimes \Lambda_{(\check{X})^{I_2}}$. Moreover the canonical maps are isomorphism 
$$\mathfrak{p}_!\Psi_{\overline{\eta}_{I_2} \to x}\mcF_{I_1 \cup I_2, W \boxtimes 1} = \mathfrak{p}_!\mcF_{I_1, W} \otimes \Lambda = \Psi_{\overline{\eta}_{I_2} \to x}\mathfrak{p}_!\mcF_{I_1 \cup I_2, W \boxtimes 1}.$$
\end{lem}

\begin{proof}
This is an immediate application of the Künneth formula for nearby cycles \ref{thmKunnethCyclesProches}. 
\end{proof}

\section{General Case}\label{sectionGeneralCase}

In this section we prove Theorem \ref{thmMain2} and Theorem \ref{thmMain1}, using the technical results established in section \ref{sectionFusion}.

\subsection{Preparations}

\subsubsection{Notations for the proof of theorem \ref{thmMain1}}
We fix $\mcA$ on $(\ChtR_{\{0\}}^{nx})_x$ and $(J,V)$ as in the introduction so that for all $(I,W)$ we have a sheaf $\mcF_{I,W}$. Lemmas \ref{lemFusionNearbyCycle2} and \ref{lemFusionNearbyCycles3} provide a modification $\widetilde{S} \to X^I$, a point $s \in \widetilde{S}$ above the diagonal point $x$ and specialization map $\overline{\eta}_{I} \to s$. We want to construct an inverse to the canonical map 
\begin{equation}\label{equationMapCan}
\can : \mathfrak{p}_!\Psi_{\overline{\eta}_{I} \to s}\mcF_{I,W} \to \Psi_{\overline{\eta}_{I} \to s}\mathfrak{p}_!\mcF_{I,W}.
\end{equation}
In section \ref{subsectionConstructionInverse}, we construct a map $\Psi_{\overline{\eta}_{I} \to s}\mathfrak{p}_!\mcF_{I,W} \to \mathfrak{p}_!\Psi_{\overline{\eta}_{I} \to s}\mcF_{I,W}$. In section \ref{subsectionZorroArgument}, we show that this map is the inverse of the map $\can$. The structure of the argument uses Zorro's lemma in a crucial way and the argument is very similar to the one appearing in \cite{XueSmoothness} and \cite[Theorem 5.2]{Salmon}. 

\subsubsection{Proof of theorem \ref{thmMain2} assuming theorem \ref{thmMain1}}

Assume that the map (\ref{equationMapCan}) is an isomorphism, and that we have chosen a total order on $I$ such that $s$ is the strict transform along the hyperplane flag.
% or the flag $(Y_i^d)$ for $d \in \N^I$ as in lemma \ref{lemStrictTransformFrob}. 
Then there is an isomorphism by construction \ref{subsubsectionIteration} and lemma \ref{lemIteratedNearbyCyclesOnshtukas}
%lemma \ref{lemFusionOnshtukas} 
of sheaves on $(\Cht_{J \cup \{0\}})_{\check{X}^J \times x}$, 
\begin{equation*}
\Psi_{\overline{\eta}_I \to s} \mcF_{I,W} \to \Psi_1\dots \Psi_n \mcF_{I,W},
\end{equation*}
where each $\Psi_i$ as in the introduction is the classical nearby cycles functor along the $i$-th projection $\Cht_{I \cup J \cup \{0\}} \to X^I \xrightarrow{\pr_i} X$. Similarly by proposition \ref{prop-p-F-I-V-constant} (1) and lemma \ref{lemFusionNearbyCycleCohomSheaves}, there is an isomorphism 
\begin{equation*}
\Psi_{\overline{\eta}_I \to s} \mathfrak{p}_!\mcF_{I,W} \to \Psi_1\dots \Psi_n \mathfrak{p}_!\mcF_{I,W},
\end{equation*}
Both isomorphisms are compatible with base change maps, that is, the following diagram is commutative:
\[\begin{tikzcd}
	{\mathfrak{p}_!\Psi_{\overline{\eta}_I \to s} \mcF_{I,W}} & {p_!\Psi_1\dots \Psi_n \mcF_{I,W}} \\
	{\Psi_{\overline{\eta}_I \to s} \mathfrak{p}_!\mcF_{I,W}} & {\Psi_1\dots \Psi_n \mathfrak{p}_!\mcF_{I,W}.}
	\arrow["\sim", from=1-1, to=1-2]
	\arrow["\can"', from=1-1, to=2-1]
	\arrow["\can", from=1-2, to=2-2]
	\arrow["\sim", from=2-1, to=2-2]
\end{tikzcd}\]
This finishes the proof of theorem \ref{thmMain2} assuming theorem \ref{thmMain1}. 

\subsection{Construction of the inverse map}\label{subsectionConstructionInverse}

We introduce three copies of $I$ which we denote by $I_1, I_2$ and $I_3$ and we denote as before by $I_{123} = I_1 \cup I_2 \cup I_3$ and by $I_{ij} = I_{i} \cup I_j$. The construction follows the example of section \ref{subsectionInversemap}. We then define the following map 

\begin{equation}\label{equationDefAlpha}
\gamma : \Psi_{\overline{\eta}_{I} \to s}\mathfrak{p}_!\mcF_{I,W} \to \mathfrak{p}_!\Psi_{\overline{\eta}_{I} \to s}\mcF_{I,W}
\end{equation}
obtained as the following composition.

\begin{align*}
\Psi_{\overline{\eta}_{I_1} \to s}\mathfrak{p}_!\mcF_{I_1,W} &= \Psi_{\overline{\eta}_{I_1} \to s}\mathfrak{p}_!\Psi_{\overline{\eta}_{I_2} \to s}\mcF_{I_{12},W \boxtimes 1} \\
&\to \Psi_{\overline{\eta}_{I_1} \to s}\mathfrak{p}_!\Psi_{\overline{\eta}_{I_2} \to s}\mcF_{I_{12},W \boxtimes (W^* \otimes W)} \\
&\to \Psi_{\overline{\eta}_{I_1} \to s}\mathfrak{p}_!\Psi_{\overline{\eta}_{I_2} \to s}\Psi_{\overline{\eta}_{I_3} \to s}\mcF_{I_{123},W \boxtimes W^* \boxtimes W} \\
&\xrightarrow{\can} \Psi_{\overline{\eta}_{I_1} \to s}\Psi_{\overline{\eta}_{I_2} \to s}\mathfrak{p}_!\Psi_{\overline{\eta}_{I_3} \to s}\mcF_{I_{123},W \boxtimes W^* \boxtimes W} \\
&\to \Psi_{\overline{\eta}_{I_2} \to s}\mathfrak{p}_!\Psi_{\overline{\eta}_{I_3} \to s}\mcF_{I_{23},(W \otimes W^*) \boxtimes W} \\
&\to \Psi_{\overline{\eta}_{I_2} \to s}\mathfrak{p}_!\Psi_{\overline{\eta}_{I_3} \to s}\mcF_{I_{23},1 \boxtimes W} \\
&= \mathfrak{p}_!\Psi_{\overline{\eta}_{I_3} \to s}\mcF_{I_{3}, W}. \\
\end{align*}

The maps are the following ones. 
\begin{enumerate}
\item The first one comes from lemma \ref{lemTrivialModification}. 
\item The second one comes from the functioriality $1 \to W^* \otimes W$ of Satake sheaves. 
\item The third one is the composition of the fusion isomorphism of nearby cycles on shtukas of proposition \ref{lemFusionNearbyCycle2} and the fusion property of the Satake sheaves, namely that $\mcF_{I_{12}, W \boxtimes (W^* \otimes W)} = (\Delta^{2, 3})^* F_{I_{123}, W \boxtimes W^* \boxtimes W}$, where $\Delta^{2, 3}$ is the partial diagonal along the last two copies of $I$. 
\item The fourth one is the canonical base change map. 
\item The fifth one is the composition of the fusion isomorphism of nearby cycles on the curve of Proposition \ref{prop-p-Psi-F-I-V-constant} (1) and lemma \ref{lemFusionNearbyCycleCohomSheaves}, and the fusion property of the Satake sheaves, namely that $\mcF_{I_{23}, (W \otimes W^*) \boxtimes 
W} = (\Delta^{1, 2})^* F_{I_{123}, W \boxtimes W^* \boxtimes W}$, where $\Delta^{1, 2}$ is the partial diagonal along the first two copies of $I$.  
\item The sixth one comes from the functoriality $W \otimes W^* \to 1$ of Satake sheave.
\item The last one is the inverse of the first one. 
\end{enumerate}

\subsection{Zorro's lemma argument}\label{subsectionZorroArgument}

Proceeding as in \cite{Salmon}, we show that the the map (\ref{equationDefAlpha}) is the inverse of $\can$. 
\begin{lem}
The composition $\gamma \circ \can$ is an isomorphism. 
\end{lem}

\begin{proof}
Consider the diagram 
\[\begin{tikzcd}
	{\mathfrak{p}_!\Psi_{\overline{\eta}_{I_1} \to s} \Psi_{\overline{\eta}_{I_2} \to s}\mcF_{I_{12},W \boxtimes 1}} & {\Psi_{\overline{\eta}_{I_1} \to s} \mathfrak{p}_!\Psi_{\overline{\eta}_{I_2} \to s}\mcF_{I_{12},W \boxtimes 1}} \\
	{\mathfrak{p}_!\Psi_{\overline{\eta}_{I_1} \to s} \Psi_{\overline{\eta}_{I_2} \to s}\mcF_{I_{12},W \boxtimes (W^* \otimes W) }} & {\Psi_{\overline{\eta}_{I_1} \to s}\mathfrak{p}_! \Psi_{\overline{\eta}_{I_2} \to s}\mcF_{I_{12},W \boxtimes (W^* \otimes W) }} \\
	{\mathfrak{p}_!\Psi_{\overline{\eta}_{I_1} \to s} \Psi_{\overline{\eta}_{I_2} \to s}\Psi_{\overline{\eta}_{I_3} \to s}\mcF_{I_{123},W \boxtimes W^* \boxtimes W }} & {\Psi_{\overline{\eta}_{I_1} \to s}\mathfrak{p}_! \Psi_{\overline{\eta}_{I_2} \to s}\Psi_{\overline{\eta}_{I_3} \to s}\mcF_{I_{123},W \boxtimes W^* \boxtimes W }} \\
	& {\Psi_{\overline{\eta}_{I_1} \to s}\Psi_{\overline{\eta}_{I_2} \to s}\mathfrak{p}_! \Psi_{\overline{\eta}_{I_3} \to s}\mcF_{I_{123},W \boxtimes W^* \boxtimes W }} \\
	{\mathfrak{p}_! \Psi_{\overline{\eta}_{I_2} \to s}\Psi_{\overline{\eta}_{I_3} \to s}\mcF_{I_{23},(W \otimes W^*) \boxtimes W }} & {\Psi_{\overline{\eta}_{I_2} \to s}\mathfrak{p}_! \Psi_{\overline{\eta}_{I_3} \to s}\mcF_{I_{23},(W \otimes W^*) \boxtimes W }} \\
	{\mathfrak{p}_! \Psi_{\overline{\eta}_{I_2} \to s}\Psi_{\overline{\eta}_{I_3} \to s}\mcF_{I_{23},1 \boxtimes W }} & {\Psi_{\overline{\eta}_{I_2} \to s}\mathfrak{p}_! \Psi_{\overline{\eta}_{I_3} \to s}\mcF_{I_{23},1 \boxtimes W }}
	\arrow["\can"{description}, from=1-1, to=1-2]
	\arrow[from=1-1, to=2-1]
	\arrow[from=1-2, to=2-2]
	\arrow["\can"{description}, from=2-1, to=2-2]
	\arrow[from=2-1, to=3-1]
	\arrow[from=2-2, to=3-2]
	\arrow["\can"{description}, from=3-1, to=3-2]
	\arrow[from=3-1, to=4-2]
	\arrow[from=3-1, to=5-1]
	\arrow["\can"{description}, from=3-2, to=4-2]
	\arrow[from=4-2, to=5-2]
	\arrow["\can"{description}, from=5-1, to=5-2]
	\arrow[from=5-1, to=6-1]
	\arrow[from=5-2, to=6-2]
	\arrow["\can"{description}, from=6-1, to=6-2]
\end{tikzcd}\]
where the right column is the composition defining $\gamma$ minus the first and last map and the left column is from top to bottom 
\begin{enumerate}
\item the map induced by the functoriality $1 \to W^* \otimes W$, 
\item the next two maps are the two vertical maps of the diagram of lemma \ref{lemFusionNearbyCycles3},
\item the last map is induced by the functoriality $W \otimes W^* \to 1$. 
\end{enumerate}
The top and bottom squares are commutative by the functoriality of the base change map for nearby cycles. The two middle squares are commutative by lemma \ref{lemFusionNearbyCycles3}. The bottom map is an isomorphism by lemma \ref{lemTrivialModification}, hence to prove the lemma it is enough to show that the composition in the left column is an isomorphism. By lemma  \ref{lemFusionNearbyCycles3}, there is a commutative diagram
\[\begin{tikzcd}
	{\mathfrak{p}_!\Psi_{\overline{\eta}_{I_1} \to s}\mcF_{I_{1},W \otimes 1}} & {\mathfrak{p}_!\Psi_{\overline{\eta}_{I_1} \to s} \Psi_{\overline{\eta}_{I_2} \to s}\mcF_{I_{12},W \boxtimes 1}} \\
	{\mathfrak{p}_!\Psi_{\overline{\eta}_{I_1} \to s}\mcF_{I_{1},W \otimes W^* \otimes W}} & {\mathfrak{p}_!\Psi_{\overline{\eta}_{I_1} \to s} \Psi_{\overline{\eta}_{I_2} \to s}\mcF_{I_{12},W \boxtimes (W^* \otimes W) }} \\
	& {\mathfrak{p}_!\Psi_{\overline{\eta}_{I_1} \to s} \Psi_{\overline{\eta}_{I_2} \to s}\Psi_{\overline{\eta}_{I_3} \to s}\mcF_{I_{123},W \boxtimes W^* \boxtimes W }} \\
	{\mathfrak{p}_!\Psi_{\overline{\eta}_{I_3} \to s}\mcF_{I_{3},V \otimes V^* \otimes V}} & {\mathfrak{p}_! \Psi_{\overline{\eta}_{I_2} \to s}\Psi_{\overline{\eta}_{I_3} \to s}\mcF_{I_{23},(W \otimes W^*) \boxtimes W }} \\
	{\mathfrak{p}_!\Psi_{\overline{\eta}_{I_3} \to s}\mcF_{I_{3},1 \otimes W}} & {\mathfrak{p}_! \Psi_{\overline{\eta}_{I_2} \to s}\Psi_{\overline{\eta}_{I_3} \to s}\mcF_{I_{23},1 \boxtimes W }}
	\arrow[from=1-1, to=1-2]
	\arrow[from=1-1, to=2-1]
	\arrow[from=1-2, to=2-2]
	\arrow[from=2-1, to=2-2]
	\arrow[from=2-1, to=3-2]
	\arrow[Rightarrow, no head, from=2-1, to=4-1]
	\arrow[from=2-2, to=3-2]
	\arrow[from=3-2, to=4-2]
	\arrow[from=4-1, to=3-2]
	\arrow[from=4-1, to=4-2]
	\arrow[from=4-1, to=5-1]
	\arrow[from=4-2, to=5-2]
	\arrow[from=5-1, to=5-2]
\end{tikzcd}\]
where all the horizontal and slanted maps are isomorphisms. The left composition in the diagram is then simply the composition induced by the functoriality of Satake sheaves 
$$V \xrightarrow{\id \otimes \coev} W \otimes W^* \otimes W \xrightarrow{\ev \otimes \id} W$$
which is an isomorphism by Zorro's lemma.   
\end{proof}

\begin{lem}
The composition $\can \circ \gamma$ is an isomorphism. 
\end{lem}

\begin{proof}
Similarly consider the following diagram 
\[\begin{tikzcd}
	{\Psi_{\overline{\eta}_{I_1} \to s}\mathfrak{p}_! \Psi_{\overline{\eta}_{I_2} \to s}\mcF_{I_{12},W \boxtimes 1}} & {\Psi_{\overline{\eta}_{I_1} \to s} \Psi_{\overline{\eta}_{I_2} \to s}\mathfrak{p}_!\mcF_{I_{12},W \boxtimes 1}} \\
	{\Psi_{\overline{\eta}_{I_1} \to s} \mathfrak{p}_!\Psi_{\overline{\eta}_{I_2} \to s}\mcF_{I_{12},W \boxtimes (W^* \otimes W) }} & {\Psi_{\overline{\eta}_{I_1} \to s} \Psi_{\overline{\eta}_{I_2} \to s}\mathfrak{p}_!\mcF_{I_{12},W \boxtimes (W^* \otimes W) }} \\
	{\Psi_{\overline{\eta}_{I_1} \to s}\mathfrak{p}_! \Psi_{\overline{\eta}_{I_2} \to s}\Psi_{\overline{\eta}_{I_3} \to s}\mcF_{I_{123},W \boxtimes W^* \boxtimes W }} \\
	{\Psi_{\overline{\eta}_{I_1} \to s} \Psi_{\overline{\eta}_{I_2} \to s}\mathfrak{p}_!\Psi_{\overline{\eta}_{I_3} \to s}\mcF_{I_{123},W \boxtimes W^* \boxtimes W }} & {\Psi_{\overline{\eta}_{I_1} \to s} \Psi_{\overline{\eta}_{I_2} \to s}\Psi_{\overline{\eta}_{I_3} \to s}\mathfrak{p}_!\mcF_{I_{123},W \boxtimes W^* \boxtimes W }} \\
	{ \Psi_{\overline{\eta}_{I_2} \to s}\mathfrak{p}_!\Psi_{\overline{\eta}_{I_3} \to s}\mcF_{I_{23},(W \otimes W^*) \boxtimes W }} & { \Psi_{\overline{\eta}_{I_2} \to s}\Psi_{\overline{\eta}_{I_3} \to s}\mathfrak{p}_!\mcF_{I_{23},(W \otimes W^*) \boxtimes W }} \\
	{ \Psi_{\overline{\eta}_{I_2} \to s}\mathfrak{p}_!\Psi_{\overline{\eta}_{I_3} \to s}\mcF_{I_{23},1 \boxtimes W }} & { \Psi_{\overline{\eta}_{I_2} \to s}\Psi_{\overline{\eta}_{I_3} \to s}\mathfrak{p}_!\mcF_{I_{23},1 \boxtimes W }}
	\arrow[from=1-1, to=1-2]
	\arrow[from=1-1, to=2-1]
	\arrow[from=1-2, to=2-2]
	\arrow[from=2-1, to=2-2]
	\arrow[from=2-1, to=3-1]
	\arrow[from=2-2, to=4-2]
	\arrow[from=3-1, to=4-1]
	\arrow[from=3-1, to=4-2]
	\arrow[from=4-1, to=4-2]
	\arrow[from=4-1, to=5-1]
	\arrow[from=4-2, to=5-2]
	\arrow[from=5-1, to=5-2]
	\arrow[from=5-1, to=6-1]
	\arrow[from=5-2, to=6-2]
	\arrow[from=6-1, to=6-2]
\end{tikzcd}\]
where the left column is the sequence of maps defining $\gamma$, the horizontal maps are the canonical base change maps and the right column is the following sequence of maps 
\begin{enumerate}
\item the first and last one come from the functorialities $1 \to W^* \otimes W$ and $W \otimes W^* \to 1$, 
\item the two middles are the fusion maps for nearby cycles proposition \ref{prop-p-F-I-V-constant} (1) and lemma \ref{lemFusionNearbyCycleCohomSheaves} with respect to $I_{23} \to I_2$ and $I_{12} \to I_2$ inducing the partial diagonal respectively. 
\end{enumerate}
The top and bottom squares of the diagram are commutative by functoriality of the base change maps, the top middle square is commutative by lemma \ref{lemFusionNearbyCycles3} and the second one by functoriality of fusion on the curve by lemma \ref{lemFusionNearbyCycleCohomSheaves}. By lemma \ref{lemTrivialModification}, the top and bottom base change maps are isomorphisms hence the lemma reduces down to showing that the composition in the right column is an isomorphism. Further applying fusion on the curve produces a commutative diagram of isomorphisms
\[\begin{tikzcd}
	{\Psi_{\overline{\eta}_{I_1} \to s} \Psi_{\overline{\eta}_{I_2} \to s}\mathfrak{p}_!\mcF_{I_{12},W \boxtimes 1}} & {\Psi_{\overline{\eta}_{I_1} \to s} \mathfrak{p}_!\mcF_{I_{1},W \otimes 1}} \\
	{\Psi_{\overline{\eta}_{I_1} \to s} \Psi_{\overline{\eta}_{I_2} \to s}\mathfrak{p}_!\mcF_{I_{12},W \boxtimes (W^* \otimes W) }} & {\Psi_{\overline{\eta}_{I_1} \to s}\mathfrak{p}_!\mcF_{I_{1},W \otimes W^* \otimes W }} \\
	{\Psi_{\overline{\eta}_{I_1} \to s} \Psi_{\overline{\eta}_{I_2} \to s}\Psi_{\overline{\eta}_{I_3} \to s}\mathfrak{p}_!\mcF_{I_{123},W \boxtimes W^* \boxtimes W }} \\
	{ \Psi_{\overline{\eta}_{I_2} \to s}\Psi_{\overline{\eta}_{I_3} \to s}\mathfrak{p}_!\mcF_{I_{23},(W \otimes W^*) \boxtimes W }} & {\Psi_{\overline{\eta}_{I_3} \to s}\mathfrak{p}_!\mcF_{I_{3},W \otimes W^* \otimes W }} \\
	{ \Psi_{\overline{\eta}_{I_2} \to s}\Psi_{\overline{\eta}_{I_3} \to s}\mathfrak{p}_!\mcF_{I_{23},1 \boxtimes W }} & {\Psi_{\overline{\eta}_{I_3} \to s}\mathfrak{p}_!\mcF_{I_{3},1 \otimes W }.}
	\arrow[from=1-1, to=1-2]
	\arrow[from=1-1, to=2-1]
	\arrow[from=1-2, to=2-2]
	\arrow[from=2-1, to=2-2]
	\arrow[from=2-1, to=3-1]
	\arrow[Rightarrow, no head, from=2-2, to=4-2]
	\arrow[from=3-1, to=2-2]
	\arrow[from=3-1, to=4-1]
	\arrow[from=3-1, to=4-2]
	\arrow[from=4-1, to=4-2]
	\arrow[from=4-2, to=5-2]
	\arrow[from=4-1, to=5-1]
	\arrow[from=5-1, to=5-2]
\end{tikzcd}\]
The right column in this last diagram is isomorphic to the identity by Zorro's lemma again. 
\end{proof}

\appendix

\section{Higher nearby cycles}\label{AppendixCyclesProches}

We recall some properties and theorems for higher nearby cycles (also called nearby cycles over a base of dimension $>1$ or nearby cycles over a general base). The idea follows from Deligne and first written by Laumon in \cite{Laumon}. Main results are developed by Orgogozo (with the help of Gabber) in \cite{Orgogozo}. These results are also explained in a survey by Illusie \cite{Illusiesurvey}. Besides, we use in a crucial way a theorem of Gabber proved by using \cite{HansenScholze} (that we will recall in Theorem \ref{thmGabber} below).

In this appendix, we fix $f : \mcY \to S$ a morphism of finite presentation between two qcqs schemes with $\ell$ invertible on $S$ and $\Lambda$ a finite ring of order a power of $\ell$. We will denote by $\mcY \overleftarrow{\times}_S S$ the corresponding oriented product \cite{TravauxDeGabber} and by $\Psi_f : \mcY \to \mcY \overleftarrow{\times}_S S$ the induced map of toposes. Similarly if $U \subset S$ is an open, we denote by $\mcY_U = \mcY \times_S U$ and by $\Psi_{U,f} : \mcY_U \to \mcY \overleftarrow{\times}_S U$ the induced map. 

\begin{rque}\label{rquePointsOfOrientedTopoi}
We recall, see \emph{loc. cit.} for a more detailled discussion, that points of the topos $\mcY \overleftarrow{\times}_S S$ are triples $(x,\phi,t)$ where $x \to \mcY$ is a geometric point whose image in $S$ we shall denote by $s$, $t$ is a geometric point of $S$ and $\phi : t \to s$ is a specialization map.
\end{rque}

Given $g : T \to S$ a map of qcqs schemes we introduce the following notations for the pullbacks to $T$,
\[\begin{tikzcd}
	{\mcY_{U_T}} & {\mcY_T} \\
	{U_T} & T \\
	& {\mcY_U} & \mcY \\
	& U & {S.}
	\arrow["{j_T}"{description}, from=1-1, to=1-2]
	\arrow["{f_{U_T}}"{description}, from=1-1, to=2-1]
	\arrow["{f_T}"{description}, from=1-2, to=2-2]
	\arrow["p"{description}, from=1-2, to=3-3]
	\arrow[from=2-1, to=2-2]
	\arrow["{g_U}"{description}, from=2-1, to=4-2]
	\arrow[from=2-2, to=4-3]
	\arrow["j"{description}, from=3-2, to=3-3]
	\arrow["{f_U}"{description}, from=3-2, to=4-2]
	\arrow["f"{description}, from=3-3, to=4-3]
	\arrow[from=4-2, to=4-3]
	\arrow["p_U", from=1-1, to=3-2]
\end{tikzcd}\]
There are two commutative diagram of toposes 
\[\begin{tikzcd}
	{\mcY_T} & {\mcY_T \overleftarrow{\times}_T T} \\
	\mcY & {\mcY \overleftarrow{\times}_S S}
	\arrow["{\Psi_{f_T}}", from=1-1, to=1-2]
	\arrow["p"', from=1-1, to=2-1]
	\arrow["{\overleftarrow{p}}", from=1-2, to=2-2]
	\arrow["{\Psi_f}"', from=2-1, to=2-2]
\end{tikzcd}\]
and 
\[\begin{tikzcd}
	{\mcY_{U_T}} & {\mcY_T \overleftarrow{\times}_T U_T} \\
	{\mcY_U} & {\mcY \overleftarrow{\times}_S U}
	\arrow["{\Psi_{U_T,f_T}}", from=1-1, to=1-2]
	\arrow["p_U"', from=1-1, to=2-1]
	\arrow["{\overleftarrow{p_U}}", from=1-2, to=2-2]
	\arrow["{\Psi_{U,f}}"', from=2-1, to=2-2]
\end{tikzcd}\]

\begin{defi}[$\Psi$-good sheaves]  \label{def-Psigood}
\begin{enumerate}
\item Let $\mcA \in \DD^b_c(\mcY, \Lambda)$, the sheaf $\mcA$ is $\Psi$-good if for all $T \to S$ the natural map 
$$\overleftarrow{p}^*\Psi_{f,*}\mcA \to \Psi_{f_T,*}p^*\mcA$$
is an isomorphism. 
\item Let $\mcA \in \DD^b_c(\mcY_U, \Lambda)$, the sheaf $\mcA$ is $(U, \Psi)$-good if for all $T \to S$ the natural map
$$\overleftarrow{p_U}^*\Psi_{U,f,*}\mcA \to \Psi_{U_T,f_T,*}p_U^*\mcA$$
is an isomorphism.
\end{enumerate}
\end{defi}

\begin{rque}
The $\Psi$-goodness hypothesis can be reformulated as 'the formation of higher nearby cycles commutes with arbitrary base change'. 
\end{rque}

\begin{thm}\cite[Théorème 1.1, Théorème 6.1]{Orgogozo}\label{thmOrgogozo}
Let $\mcA \in \DD_c^b(\mcY, \Lambda)$, there exists a modification $g : S' \to S$ such that $p^*\mcA$ is $\Psi$-good and $\Psi_{f_{S'},*}\mcA$ is constructible in the sense of \emph{loc. cit.}
\end{thm}

\begin{rque}
In particular if $S$ is the spectrum of a valuation ring, then any sheaf over it is $\Psi$-good as any modification of valuative schemes have sections.
\end{rque} 

\begin{corol}\label{corolOrgogozo}
Let $\mcA \in \DD^b_c(\mcY_U, \Lambda)$, there exists a modification $g : S' \to S$ such that $p_U^*\mcA$ is $(U, \Psi)$-good and $\Psi_{U_{S'},f_{S'},*}\mcA$ is constructible in the sense of \emph{loc. cit.}
\end{corol}

\begin{rque} \label{rqueOrgogozo}
Consider the map of toposes $h : \mcY \overleftarrow{\times}_S S \to S \overleftarrow{\times}_S S$ induced by the map $\mcY \to S$, according to remark \ref{rquePointsOfOrientedTopoi}, points of $ S \overleftarrow{\times}_S S$ are exactly specialization maps of geometric points of $S$. Let $(t \to s) \in  S \overleftarrow{\times}_S S$ be such a point, the fiber over the point $t \to s$ is the subtopos of $\mcY \overleftarrow{\times}_S S$ whose points are geometric points $\bar{x} \to \mcY$ such that the image of $\bar{x}$ in $S$ is $s$. There is an identification $h^{-1}(t \to s) \simeq \mcY_s$. 

Let us now make the relation between $\Psi_{f,*}$ and $\Psi_{t \to s}$ (defined in \ref{subsectionReminderNearbyCycles}
). Let $\mcA \in \DD^b_c(\mcY, \Lambda)$ be $\Psi$-good with respect to $S$. Then by \cite[Théorème 5.1]{Orgogozo} the canonical morphism of sheaves on $\mcY_s$ 
$$(\Psi_{f,*}\mcA)_{|h^{-1}(t \to s)} \to \Psi_{t \to s}\mcA$$ 
is an isomorphism. The $\Psi$-goodness is necessary for this assumption as the LHS in the above equation computes the cohomology of the Milnor tubes while the RHS computes the cohomology of the Milnor fibers, we refer to \emph{loc. cit.} for a discussion.
\end{rque}

\begin{rque}  \label{rquebasechange}
Assume that $\mcA \in \DD^b_c(\mcY, \Lambda)$ is $\Psi$-good. Then for all $T \to S$ and all specialization maps $\alpha_T : a_T \to b_T$ in $T$ with image $\alpha :a \to b$ in $S$, denote by $p_b : (\mcY_T)_{b_T} \to \mcY_b$ then the natural base change map
$$p_b^*\Psi_{a \to b}\mcA \to \Psi_{a_T \to b_T}p^*\mcA$$
is an isomorphism. 

Assume that $\mcA \in \DD^b_c(\mcY_U, \Lambda)$ is $(U,\Psi)$-good, then the same statement holds whenever $a_T \in U_T$. 
\end{rque}

\begin{proof}[Proof of corollary \ref{corolOrgogozo}]
After replacing $S$ by a modification $S'$ such that $j_!\mcA$ is $\Psi$-good, where $j : \mcY_U \to \mcY$ is the open inclusion, we can assume that $j_!\mcA$ is $\Psi$-good. Let $T \to S$, it is enough to check that for all specializations $a_T \to b_T$ in $T$ with $a_T \in U_T$ the induced map 
$$p_b^*\Psi_{a \to b}\mcA \to \Psi_{a_T \to b_T}p^*\mcA$$
is an isomorphism. This is immediate since $j_!\mcA$ is $\Psi$-good. 
\end{proof}

Let $s \xrightarrow{\alpha} t \xrightarrow{\beta} u$ be a sequence of specialization maps in $S$ and consider the induced diagram 
\[\begin{tikzcd}
	{\mcY_s} & {\mcY_t} & {\mcY_u} \\
	{\mcY_{S_{(s)}}} & {\mcY_{S_{(t)}}} & {\mcY_{S_{(u)}}}
	\arrow["{i_s}"{description}, from=1-1, to=2-1]
	\arrow["{i_t}"{description}, from=1-2, to=2-2]
	\arrow["{i_u}"{description}, from=1-3, to=2-3]
	\arrow["{j_{\alpha}}"{description}, from=2-1, to=2-2]
	\arrow["{j_{\beta}}"{description}, from=2-2, to=2-3]
\end{tikzcd}\]
The adjunction map $\id \to i_{t,*}i_t^*$ then induces a canonical map 
$$\Psi_{s \to u} = i_u^*j_{\beta,*}j_{\alpha,*}i_{s,*} \to i_u^*j_{\beta,*}i_{t,*}i_t^*j_{\alpha,*}i_{s,*} = \Psi_{t \to u}\Psi_{s \to t}.$$

\begin{thm}[Gabber]\label{thmGabber}
Let $\mcA \in \DD_c^b(\mcY, \Lambda)$. There exists a modification $S' \rightarrow S$ such that for all specialization maps $s \rightarrow t \rightarrow u$ in $S'$, the natural map 
\begin{equation}\label{equationthmGabber}
\Psi_{s \rightarrow u}\mcA \rightarrow \Psi_{t \rightarrow u} \Psi_{s \rightarrow t} \mcA
\end{equation}
is an isomorphism. 
\end{thm}

\begin{proof}
This theorem is stated in \protect{\cite[Theorem 4.5]{Abe}} without proof.
Since we did not find a proof in the literature we provide one here. By theorem \ref{thmOrgogozo}, we can modify $S$ and therefore assume that $\mcA$ is $\Psi$-good. Let $s \xrightarrow{\alpha} t \xrightarrow{\beta} u$ be a chain of specializations in $S$. There exists a rank $2$ absolutely integrally closed valuation ring $V$ and a map $h : \Spec(V) \to S$ representing this chain of specializations. Since $\mcA$ is $\Psi$-good, we can replace $S$ by $\Spec(V)$. Let $\eta \in \Spec(V)$ be the generic point and $j : \eta \to \Spec(V)$ and both sides of (\ref{equationthmGabber}) depend only on $j^*\mcA$. By \cite[Theorem 4.1]{HansenScholze}, the sheaf $j_*j^*\mcA$ is ULA relative to $\Spec(V)$ and for ULA sheaves both sides are canonically isomorphic to $(j_*j^*\mcA)_u$ hence the map (\ref{equationthmGabber}) is an isomorphism.
\end{proof}

\begin{rque}
It should be noted from the proof that the same modification provided by theorem \ref{thmOrgogozo} works for theorem \ref{thmGabber}. 
\end{rque} 

\begin{thm}[Künneth formula, \protect{\cite[Theorem 2.3]{Illusie}}]\label{thmKunnethCyclesProches}
Assume that $S = S_1 \times S_2$, that $\mcY = \mcY_1 \times \mcY_2$ and $f = f_1 \times f_2$ for two maps $f_1 : \mcY_1 \to S_1$ and $f_2 : \mcY_2 \to S_2$. Then for $\mcA_i \in \DD^b_c(\mcY_i, \Lambda)$, there is a canonical map 
$$\overleftarrow{\Psi}_{f_1,*}\mcA_1 \boxtimes \overleftarrow{\Psi}_{f_2,*}\mcA_2 \to \overleftarrow{\Psi}_{f,*}(\mcA_1 \boxtimes \mcA_2).$$
Furthermore, if both  $\mcA_1$ and $\mcA_2$ are $\Psi$-good relative to $S_1$ and $S_2$ respectively then $\mcA_1 \boxtimes \mcA_2$ is $\Psi$-good with respect to $S$ and the above map is an isomorphism. 
\end{thm}

\bibliographystyle{alpha}
\bibliography{biblio.bib}

\end{document}